\documentclass[12pt]{article}
\usepackage{amssymb}
\usepackage{latexsym}

\newtheorem{theorem}{Theorem}[section]
\newtheorem{definition}[theorem]{Definition}
\newtheorem{proposition}[theorem]{Proposition}
\newtheorem{lemma}[theorem]{Lemma}

\newtheorem{remark}[theorem]{Remark}

\newtheorem{claim}{Claim}
\newtheorem{question}[theorem]{Question}
\newcommand {\z}{{\mathbb Z}}

\newcommand {\mod}{\mbox{mod }}

\renewcommand {\thefootnote}{\fnsymbol{footnote}}

\begin{document}

\baselineskip=16pt

\begin{center}
{\bf\Large
Detailed Structure for\\
\vspace{0.1in}
Freiman's $3k-3$ Theorem}
\end{center}

\bigskip

\begin{center}
Renling Jin\\
\vspace{0.1in}
College of Charleston\\
{\tt jinr@cofc.edu}\\
July, 2013
\end{center}

\bigskip\bigskip

\renewcommand{\thefootnote}{\arabic{footnote}}\addtocounter{footnote}{-1}

\begin{quote}

 \centerline{Abstract}

Let $A$ be a finite set of integers. We prove that
if $|A|\geqslant 2$ and $|A+A|=3|A|-3$, then one of the following is true:
\begin{enumerate}
\item $A$ is a bi-arithmetic progression;
\item $A+A$ contains an arithmetic progression of length $2|A|-1$;
\item $|A|=6$ and $A$ is Freiman isomorphic to the set

\qquad $\{(0,0),(0,1),(0,2),(1,0),(1,1),(2,0)\}\subseteq\z^2$;
\item $|A|=1+(\max A-\min A)/(2\gcd(A-\min A))$ and
$A$ is Freiman isomorphic to a set in either the form of

\qquad $\{0,2,\ldots,2k\}\cup B\cup\{n\}$ \newline
for some non-negative integer $k\leqslant \frac{1}{2}n-2$ or the form of

\qquad $\{0\}\cup C\cup D\cup\{n\}$\newline
where $n=2|A|-2$, $B$ is left dense in $[2k,n-1]$, $C$ is right dense in
$[1,u]$ for some $u\in [4,n-6]$, $D$ is left dense in $[u+2,n-1]$, $B,C,D$ are
anti-symmetric and additively minimal in the correspondent
host intervals\footnote{See Definition \ref{triangle}.}.
\end{enumerate}
\end{quote}

\section{Introduction and Propositions}

Let $A$ be a finite set of integers and $|A|$ be the cardinality of $A$.
Let $A\pm B:=\{a\pm b:a\in A,\,b\in B\}$ for any sets of integers $A$ and $B$.

Freiman's inverse problem for small doubling constants
seeks structural information for $A$ or $\,2A:=A+A$ when the size of $2A$
is small, say for example, less than $4|A|$.

Let $X$ be a subset of an abelian semigroup $G$ and $Y$ be a subset of an abelian
semigroup $G'$. A bijection $\varphi:X\mapsto Y$ is called a {\it Freiman isomorphism}
of order $2$ if for any $a,b,c,d\in X$,
\[a+b=c+d\,\,\mbox{ if and only if }\,\,\varphi(a)+\varphi(b)=\varphi(c)+\varphi(d).\]
We will use only order $2$ Freiman isomorphism in this paper. Therefore,
the word ``order $2$'' will be omitted.

Let's call a set $B$ a {\it bi-arithmetic progression}
if $B$ is the Freiman isomorphism image of the set
\[B'=\{(0,0),(0,1),\ldots,(0,n_0)\}\cup\{(1,0),(1,1),\ldots,(1,n_1)\}\]
as a subset of the usual additive group $(\mathbb{Z}^2,+)$
where $n_0,n_1\geqslant 0$ and $n_0+n_1+2=|B|$
\footnote{$B=I_0\cup I_1\subseteq\mathbb{Z}$ is a bi-arithmetic progression
of difference $d$ if and only if $I_0$ and $I_1$
are two arithmetic progressions of difference $d$ and
$2I_0$, $I_0+I_1$, $2I_1$ are pairwise disjoint.}. Let $\varphi$ be the Freiman isomorphism
from $B'$ to $B$ and $I_i=\varphi(\{(i,0),(i,1),\ldots,(i,n_i)\})$
for $i=0,1$. The common difference of
$I_0$ and $I_1$ is called the difference of bi-arithmetic progression $B$.
We often call the expression $I_0\cup I_1$ a (bi-arithmetic progression) decomposition
of $B$. For example, $B=\{0,3,5,6,8\}$ is a bi-arithmetic progression of difference $3$ and
has a decomposition $\{0,3,6\}\cup\{5,8\}$.

Around 1960 Freiman proved the following two
theorems in \cite[page 11, page 15]{freiman1} or \cite{nathanson}.

\begin{theorem}[G.\ A.\ Freiman]\label{2k-1+b}

Let $A$ be a finite set of integers and $|A|>2$. If
$|2A|=2|A|-1+b<3|A|-3$, then $A$ is a subset of an arithmetic progression of length
at most $|A|+b$.

\end{theorem}

\begin{theorem}[G.\ A.\ Freiman]\label{3k-3}

Let $A$ be a finite set of integers and $|A|\geqslant 2$. If
$|2A|=2|A|-1+b=3|A|-3$, then one of the following is true.
\begin{enumerate}
\item $A$ is a bi-arithmetic progression;
\item $A$ is a subset of an arithmetic progression of length
at most $2|A|-1$;
\item $|A|=6$ and $A$ is a Freiman isomorphism image
of the set $K_6$ where
\begin{equation}\label{sixpoints}
K_6=\{(0,0),(0,1),(0,2),(1,0),(1,1),(2,0)\}\subseteq\z^2.
\end{equation}
\end{enumerate}
\end{theorem}

The conclusion that $A$ is a subset of an arithmetic progression $I$ of length
$2|A|-1$ in Theorem \ref{3k-3} indicates $|A|>\frac{1}{2}|I|$, or we can say that
$A$ is a large subset of the arithmetic progression $I$.

Part 1 and 2 in Theorem \ref{3k-3} show the regularity of the structure of $A$
 when $|2A|=3|A|-3$. We view part 3 as an exception.
If $A$ is the set $\{0,a,2a,b,b+a,2b\}$ for any $b>4a$,
then $A$ is Freiman isomorphic to $K_6$. Clearly, this $A$ can be neither an
arithmetic progression nor a bi-arithmetic progression of reasonable length
while $a$ and $b$ can be as large as we want.

Each element in $V=\{(0,2), (2,0), (0,0)\}$ is called a vertex of $K_6$.
Notice that each permutation of $V$ can be extended to a Freiman isomorphism
from $K_6$ to $K_6$. If $\varphi:K_6\mapsto B$ is a Freiman isomorphism,
we also call the elements in $\varphi(V)$ vertices of $\varphi(K_6)$.

Theorem \ref{3k-3} is much more difficult to prove than Theorem \ref{2k-1+b} is.
There has been a few generalizations of Theorem \ref{3k-3}. In \cite{HP} it is proved that
the structure of $A$ is the same as the structure of $A$
characterized in Theorem \ref{3k-3} when $|A\pm A|=3|A|-3$.
In \cite{jin2}, a generalization of Theorem \ref{3k-3} is given, which
 characterizes the structure of $A$ when $|A|$ is sufficiently large
and $|2A|=3|A|-3+b$ for $0\leqslant b\leqslant \epsilon|A|$, where $\epsilon$ is a
small positive real number.

 Recently, Freiman discovered in \cite{freiman2,freiman3} some interesting detailed structural
 information of $A$ when $|2A|<3|A|-3$.
 By saying ``detailed information'' we mean any structural information other than
that of $A$ being a large subset of an arithmetic progression.
The key result in \cite{freiman2,freiman3} is the following.

\begin{theorem}[A.~G.~Freiman, 2009]\label{freimanthm}
Let $A$ be a finite set of integers. If
\[|2A|<3|A|-3,\]
then $2A$ contains an arithmetic progression of length $2|A|-1$.
\end{theorem}

In \cite{BG} Theorem \ref{freimanthm} is generalized to the sum of two distinct sets.
Similar to that Theorem \ref{freimanthm} adds extra detailed structural information
to the structural information obtained in Theorem \ref{2k-1+b},
Theorem \ref{BGthm} adds extra detailed structural information to
the structural information obtained in \cite{LS} and \cite{stanchescu1}
for the addition of two distinct sets.

\begin{theorem}[I.~Bardaji and D.~J.~Grynkiewicz, 2010]\label{BGthm}
Let $A$ and $B$ be finite nonempty sets of integers with \[\max B-\min B\leqslant
\max A-\min A\leqslant |A|+|B|-3\]
\[\mbox{and }\,|A+B|\leqslant |A|+2|B|-3-\delta(A,B).\]
Then $A+B$ contains an interval of integers of length $|A|+|B|-1$.
\end{theorem}
The number $\delta(A,B)$ in Theorem \ref{BGthm} is defined to be $1$ if
$A+t\subseteq B$ for some integer $t$ and to be $0$ otherwise. If checking carefully, the reader
can find that the condition $\max A-\min A\leqslant |A|+|B|-3$ in
Theorem \ref{BGthm} can be weakened to $\max A-\min A\leqslant |A|+|B|-2$
when $\max B-\min B<\max A-\min A$.

\medskip

In this paper we seek {\em detailed} structural information for $A$ when $|2A|=3|A|-3$.
The information we have found is consistent with that in
Theorem \ref{freimanthm} and Theorem \ref{BGthm} with some exceptions.

Let $a$ and $b$ be integers. Throughout this paper we will write $[a,b]$ for the interval
of {\it integers} between $a$ and $b$ including $a$ and $b$. Notice that $[a,b]=\emptyset$
if $a>b$. For any set $A$ of integers,
we will use the following notation:
\[A(a,b):=|A\cap [a,b]|.\]
When $b$ is an integer, we write $b\pm A$ for $\{b\}\pm A$ and $A\pm b$ for $A\pm\{b\}$.

We now introduce a few  propositions, which will be used in the proof of the main result.

\begin{proposition}\label{pigeonhole}
If $A(x,y)>\frac{1}{2}(y-x+1)$, then $y+x\in 2A$.
\end{proposition}

\noindent {\bf Proof}\quad If $A(x,y)>\frac{1}{2}(y-x+1)$, then
$A\cap (x+y-A)\cap [x,y]\not=\emptyset$.
Let $a\in A\cap (x+y-A)\cap [x,y]$ and $a'\in A$ be such that
$a=x+y-a'$. Then $x+y=a+a'\in 2A$.

\begin{proposition}\label{sixpoints2}
If $\varphi:K_6\mapsto B$ is an Freiman isomorphism from $K_6$ in (\ref{sixpoints}) to $B$,
then
\begin{enumerate}
\item $\min B$ and $\max B$ are vertices of $\varphi(K_6)$.
\item If $x,y\in B$ are vertices, then $\frac{1}{2}(x+y)\in B$.
\item If $B\subseteq [a,b]$, then $b-a\geqslant 10$.
\item If $B\subseteq [0,10]$, then $B$ is either $B_1=\{0,1,2,5,6,10\}$,\newline
or $B_2=\{0,2,4,5,7,10\}$, or $B_3=10-B_1$, or $B_4=10-B_2$.
\end{enumerate}
\end{proposition}

\noindent {\bf Proof}\quad Part 1: If, for example, $\varphi((0,1))=\min B$, then
$\min B$ is the middle term of a three-term arithmetic progression $\varphi((0,0))$,
$\varphi((0,1))$, and $\varphi((0,2))$. Hence one of
$\varphi((0,2)),\varphi((0,0))$ must be smaller than $\min B$, which contradicts the
minimality of $\min B$. The argument for $\max B$ is similar.

Part 2 follows from the definition of Freiman isomorphism.

Part 3: Suppose, for example, $\varphi(\{(0,0),(0,1),(0,2)\})=\{a,a+d,a+2d\}$ where
$a=\min B$.
Then $d=1$ implies that $\varphi((1,0))\geqslant a+5$ and $b\geqslant\varphi((2,0))\geqslant
a+10$. Hence $b-a\geqslant 10$. If $d=2$, then $\varphi((1,0))\geqslant a+5$ and
$b\geqslant\varphi((2,0))\geqslant a+10$. Hence $b-a\geqslant 10$.
Suppose that $d=3$. If $\varphi((1,0))=a+1$, then $\varphi((2,0))=a+2$ and
$\varphi((1,1))=a+4$. But now $\varphi((1,0))+\varphi((0,1))=2\varphi((2,0))$, a contradiction
to that $\varphi$ is a Freiman isomorphism.
If $\varphi((1,0))=a+2$, then $\varphi((2,0))=a+4$ and $\varphi((1,1))
=a+5$. But now $\varphi((1,0))+\varphi((2,0))=2\varphi((0,1))$, again a contradiction
to that $\varphi$ is a Freiman isomorphism.
If $\varphi((1,0))=a+4$, then $\varphi((2,0))=a+8$ and $\varphi((1,1))=a+7$.
But now $\varphi((1,0))+\varphi((2,0))=2\varphi((0,2))$, again a contradiction.
So we can assume that $\varphi((1,0))\geqslant a+5$, which implies that $b\geqslant
\varphi((2,0))\geqslant a+10$. If $d\geqslant 4$, then $b\geqslant\varphi((0,2))=a+12$.

(4) Left to the reader to check.

\medskip

We introduce new names of some set configurations in order to be more
efficient and informative when describing
the set configuration in part 4 of Theorem \ref{main1}.

\begin{definition}\label{triangle}
Let $B\subseteq [u,v]$.
\begin{itemize}
\item $B$ is {\bf half dense} in $[u,v]$ if $B(u,v)=\frac{1}{2}(v-u+1)$;
\item $B$ is {\bf anti-symmetric} in $[u,v]$ if $B\cap(u+v-B)=\emptyset$
and $B\cup(u+v-B)=[u,v]$;
\item A half dense set $B$ in $[u,v]$ is {\bf left dense} in $[u,v]$
if $B(u,x)>\frac{1}{2}(x-u+1)$ for any $x\in [u,v-1]$;
A half dense set $B$ in $[u,v]$ is {\bf right dense} in $[u,v]$
if $B(x,v)>\frac{1}{2}(v-x+1)$ for any $x\in [u+1,v]$;
\item For a left dense and anti-symmetric set $B$ in $[u,v]$,
$B$ is {\bf additively minimal} in $[u,v]$ if \[2(B\cup\{v+1\})=[2u,u+v-1]\cup(v+1+(B\cup\{v+1\})).\]
For a right dense and anti-symmetric set $B$ in $[u,v]$,
$B$ is {\bf additively minimal} in $[u,v]$ if \[2(B\cup\{u-1\})=[u+v+1,2v]\cup(u-1+(B\cup\{u-1\})).\]
\end{itemize}
\end{definition}

We call the interval $[u,v]$ in Definition \ref{triangle}
the {\it host interval} of $B$ even though $u$ or $v$ may not be in $B$. Here
$|B|=\frac{1}{2}|[u,v]|$ is the reason why we give the name ``host interval.''

The following are some straightforward consequences of Definition \ref{triangle}.
An argument for a left dense set in the proposition below
works also for a right dense set by symmetry.

\begin{proposition}\label{addminim}
Let $B\subseteq [u,v]$.
\begin{enumerate}
\item $B$ is anti-symmetric in $[u,v]$ if and only if
$B(u,v)$ is half dense in $[u,v]$ and $u+v\not\in 2B$.
\item If $B$ is left dense in $[u,v]$ and $v-u>1$, then $u,u+1\in B$ and $v,v-1\not\in B$.
\item If $B$ is left dense in $[u,v]$, then
\[2(B\cup\{v+1\})\supseteq [2u,u+v-1]\cup(v+1+(B\cup\{v+1\})).\]
\item If $B$ is left dense in $[u,v]$, then
$|2(B\cup\{v+1\})|=3|B\cup\{v+1\}|-3$ if and only if
$B$ is anti-symmetric and additively minimal in $[u,v]$.
\item If $B$ is left dense in $[u,v]$,
then either $B$ is the interval $\left[u,\frac{1}{2}(u+v-1)\right]$ or
$|2(B\cup\{b\})|>3|B\cup\{b\}|-3$ for any $b>v+1$.
\item If $C$ is right dense, anti-symmetric, and additively minimal in $[1,u]$,
and $D$ is left dense, anti-symmetric, and additively minimal in $[u+2,n-1]$
for some $u\in [4,n-6]$, then
\[|2(\{0\}\cup C\cup D\cup\{n\})|=3|\{0\}\cup C\cup D\cup\{n\}|-3.\]
\end{enumerate}
\end{proposition}

\noindent {\bf Proof}\quad Part 1 and 2 are easy consequences of the definition.
Part 3 is a consequence of Proposition \ref{pigeonhole}.
Part 4 follows from part 3 because the cardinality of
$[2u,u+v-1]\cup(v+1+(B\cup\{v+1\}))$ is $3|B\cup\{v+1\}|-3$.

For part 5: For convenience we can assume, without loss of generality, that $u=0$.
If $B$ is an interval, then $B=\left[0,\frac{1}{2}(v-1)\right]$ because
$B$ is half dense in $[0,v]$. Suppose that $B$ is not an interval.

Let $a=\max B$. Since $B$ is not an interval, we have that
 $a>\frac{1}{2}(v-1)$, which implies $a\geqslant v-a$. If
$v-a\not\in B$, then $a>v-a$ because $a\in B$. Let $b>v+1$. If $v\in 2B$, then
\begin{eqnarray*}
\lefteqn{|2(B\cup\{b\})|\geqslant |[0,v]\cup(b+(B\cup\{b\}))|}\\
& &=2|B|+|B\cup\{b\}|=3|B|+1=3|B\cup\{b\}|-2.
\end{eqnarray*}
So we can assume that $v\not\in 2B$. This implies that
$B$ is anti-symmetric in $[0,v]$. We now want to show that $v+1\in 2B$.

Since $v\not\in 2B$, then $v-a\not\in B$. Hence
\begin{eqnarray*}
\lefteqn{\frac{1}{2}(v+1)=B(0,v)}\\
& &=B(0,v-a-1)+B(v-a+1,a)\\
& &\leqslant v-a+B(v-a+1,a),
\end{eqnarray*}
which implies that
$B(v-a+1,a)\geqslant\frac{1}{2}(2a-v+1)>\frac{1}{2}(2a-v)$. Hence
\[v+1=a+(v-a+1)\in 2(B\cap[v-a+1,a])\subseteq 2B\]
by Proposition \ref{pigeonhole}.
This again implies that $|2(B\cup\{b\})|>3A(0,b)-3$.

For part 6: Let $A=\{0\}\cup C\cup D\cup\{n\}$. Then $|A|=|C|+|D|+2
=\frac{1}{2}u+\frac{1}{2}(n-u-2)+2=\frac{1}{2}(n+2)$ and

\[|2A|=|0+(\{0\}\cup C)|+|[u+2,u+2+n-2]|+|n+(D\cup\{n\})|\]
\[=|A|+n-1=3|A|-3.\]

\begin{remark}
Part 5 of Proposition \ref{addminim} justifies the definition of a left dense
set to be additively minimal by looking at the cardinality of
$2(B\cup\{v+1\})$ instead of the cardinality of
$2(B\cup\{b\})$ for any $b>v+1$. Since $v+1$ is implicitly determined
by the definition of additive minimality, we can call $B$
additively minimal in its host interval $[u,v]$ without mentioning
the element $v+1$. Part 1 -- 5 in Proposition \ref{addminim} can be
restated for right dense cases.
\end{remark}

\noindent {\bf Blank Assumption}\quad After normalization,
we can always assume, throughout this paper,
that the set $A$ satisfies
\begin{equation}\label{condition}
0=\min A,\,\gcd(A)=1,\mbox{ and }\,n=\max A.
\end{equation}

\medskip

\begin{proposition}\label{sub3k-3prop}
Suppose that $0<a<b<n$ and $A\subseteq [0,n]$ such that $A\cap[a,b]=\{a,b\}$.
Then
\begin{enumerate}
\item Clearly,
\begin{equation}\label{intersection}
(2(A\cap[0,b]))\cap(2(A\cap[a,n]))=\{2a,a+b,2b\}.
\end{equation}
\item If $|2A|=3|A|-3$,
\[|2(A\cap[0,b])|\geqslant
3A(0,b)-3,\,\mbox{ and }\,|2(A\cap[a,n])|\geqslant 3A(a,n)-3,\]
\begin{equation}\label{sub3k-3}
\mbox{then }\,|2(A\cap[0,b])|=3A(0,b)-3,\,\,|2(A\cap[a,n])|=3A(a,n)-3,
\end{equation}
\begin{equation}\label{twopieces}
\mbox{and }\,(2A)\smallsetminus ((2(A\cap[0,b]))\cup(2(A\cap[a,n])))
\end{equation}
is an empty set.
\item Let $B\subseteq [u,v]$, $u,v\in B$, and $\gcd(B-u)=1$.
If $|B|\leqslant\frac{1}{2}(v-u+3)$, then $|2B|\geqslant 3|B|-3$. If
$|B|\leqslant\frac{1}{2}(v-u+1)$ and $|2B|= 3|B|-3$,
then $B$ is either a bi-arithmetic progression or a Freiman isomorphism image of $K_6$.
\end{enumerate}
\end{proposition}

\noindent {\bf Proof}\quad
Part 2 follows the inequalities
\begin{eqnarray*}
\lefteqn{3|A|-3=|2A|}\\
& &\geqslant |2(A\cap[0,b])|+|2(A\cap[a,n])|\\
& &\quad -|\{2a,a+b,2b\}|\\
& &\geqslant 3A(0,b)-3+3A(a,n)-3-3\\
& &=3A(0,a-1)+3+3A(a,n)-6=3|A|-3,
\end{eqnarray*}
which imply (\ref{sub3k-3}) and $|2A|=|2(A\cap[0,b])|+|2(A\cap[a,n])|-3$.

The first ``if'' sentence in part 3 follows from Theorem \ref{2k-1+b}.
The second ``if'' sentence
in part 3 simply eliminates the possibility of part 2 of Theorem \ref{3k-3}.

\section{Main Theorem}

Throughout this section, the letter $A$ always represents a finite set of integers
and {\bf satisfies (\ref{condition})}.
The following is the main theorem in this paper.

\begin{theorem}\label{main1} If $|A|\geqslant 2$ and
\begin{equation}\label{2k-1}
|2A|=3|A|-3,
\end{equation}
then one of the following must be true.
\begin{enumerate}
\item $A$ is a bi-arithmetic progression;
\item $2A$ contains an interval of length $2|A|-1$;
\item $|A|=6$ and $A$ is a Freiman isomorphism image of the set $K_6$
defined in (\ref{sixpoints}).
\item $|A|=\frac{1}{2}(n+2)$ and either
$A$ is in the form of $T_{k,n}$ or in the form of $n-T_{k,n}$, or
in the form of $S_{u,n}$ where
\begin{equation}\label{tk}
T_{k,n}=\{0,2,\ldots,2k\}\cup B\cup\{n\}\subseteq\mathbb{Z}
\end{equation}
for some $k\in \left[0,\frac{1}{2}n-2\right]$ such that
$B$ is left dense, anti-symmetric, and additively minimal in $[2k,n-1]$, and
\begin{equation}\label{su}
S_{u,n}=\{0\}\cup C\cup D\cup\{n\}
\end{equation}
for some $u\in [4,n-6]$ such that $C$ is right dense,
anti-symmetric, and additively minimal in $[1,u]$, and
$D$ is left dense, anti-symmetric, and additively minimal in $[u+2,n-1]$.

\end{enumerate}
\end{theorem}

\begin{remark}
Notice that $T_{k,n}$ and $S_{u,n}$ are not unique sets but sets in collections.
For example, $\{0,1,2,3,4,10\}$, $\{0,1,3,4,7,10\}$,
$\{0,1,2,5,6,10\}$, $\{0,1,2,4,6,10\}$ are all
in the form of $T_{0,10}$.

Notice also that
$2A$ for the set $A$ in part 4 of Theorem \ref{main1}
contains an interval of length $2|A|-3$ because $2T_{k,n}$ contains the interval
$[2k,2k+n-2]$ and $2S_{u,n}$ contains the interval $[u+2,n+u]$.

The structures described in four parts of Theorem \ref{main1} are not mutually exclusive.
\end{remark}

\medskip

\noindent {\bf Proof of Theorem \ref{main1}}:\quad
Without loss of generality, we can assume that
$A$ is not a bi-arithmetic progression and is not a Freiman isomorphism image
of $K_6$ in (\ref{sixpoints}). By Theorem \ref{3k-3},
we have that $n+1\leqslant 2|A|-1$ or equivalently, $|A|\geqslant\frac{1}{2}(n+2)$.
It suffices to show that $A$ satisfies part 2 or part 4 of Theorem \ref{main1}.

Let $H=[0,n]\smallsetminus A$ and $h=|H|$.
The elements in $H$ are called the holes of $A$.
Thus $h$ counts the number of holes in $H$. A non-empty
interval $[x,y]\subseteq H$ is called a gap of $A$ if $x-1,y+1\in A$. We now divide
the proof into two parts and devote one subsection for each part.

\subsection{Proof of Theorem \ref{main1} when $|A|>\frac{1}{2}(n+2)$}

In this subsection we show that part 2 of Theorem \ref{main1} is true.
Notice that $A$ in part 4 has cardinality $\frac{1}{2}(n+2)$.
Hence part 4 is irrelevant in this subsection.

For each $x\in [0,n]$, it is true that
\begin{equation}\label{bothsides}
\mbox{either }\,A(0,x)>\frac{1}{2}(x+1)\,\mbox{ or }\,
A(x,n)>\frac{1}{2}(n-x+1)
\end{equation}
 because otherwise
\[|A|\leqslant
A(0,x)+A(x,n)\leqslant\frac{1}{2}(x+1)+\frac{1}{2}(n-x+1)\leqslant\frac{1}{2}n+1,\]
which contradicts the assumption that $|A|>\frac{1}{2}n+1$. So
for any $x\in [0,n]$, either $x\in 2A$ or
$x+n\in 2A$ by Proposition \ref{pigeonhole}. Let
\begin{itemize}
\item $H_1=\{x\in H:x\not\in 2A\,\mbox{ and }\,x+n\in 2A\}$ and $h_1=|H_1|$,
\item $H_2=\{x\in H:x\in 2A\,\mbox{ and }\,x+n\not\in 2A\}$ and
$h_2=|H_2|$,
\item $H_3=\{x\in H:x\in 2A\,\mbox{ and }\,x+n\in 2A\}$
and $h_3=|H_3|$.
\end{itemize}
In \cite{freiman2}, the elements in $H_1$ are called left stable holes,
the elements in $H_2$ are called
right stable holes, and the elements in $H_3$ are called unstable holes.
By (\ref{bothsides}) we have that $H=H_1\cup H_2\cup H_3$ and $h=h_1+h_2+h_3$.

Since $|A\cup (n+A)|=2|A|-1$, we have that \[|A|-2=|(2A)\smallsetminus (A\cup (n+A))|\]
by (\ref{2k-1}). It is easy to verify that three sets
$B_1=\{x+n:x\in H_1\}$, $B_2=\{x:x\in H_2\}$, and
$B_3=\{x,x+n:x\in H_3\}$ are pairwise
disjoint and $B_1\cup B_2\cup B_3=(2A)\smallsetminus (A\cup (n+A))$.
Hence $|A|-2=h_1+h_2+2h_3$, which implies that
\begin{equation}\label{wholeh3}
|A|-2-h=|A|-2-h_1-h_2-h_3=h_3.
\end{equation}

We now prove the following lemma which implies that $2A$ contains $2|A|-1$
consecutive integers.

\begin{lemma}\label{lr}
Let $l,r\in [0,n]$ be such that $A(0,l)\leqslant\frac{1}{2}(l+1)$
and $A(n-r,n)\leqslant\frac{1}{2}(r+1)$. Then $l<n-r$.
\end{lemma}

Lemma \ref{lr} implies part 2 of Theorem \ref{main1} by the following argument. Let $l'$ be the
greatest integer in $[-1,n]$ such that $A(0,l')\leqslant\frac{1}{2}(l'+1)$ and
$r'$ be the greatest integer in $[-1,n]$ such that $A(n-r',n)\leqslant\frac{1}{2}(r'+1)$.
Then $l'<n-r'$ by Lemma \ref{lr}. For each $x>l'$, it is true that
$A(0,x)>\frac{1}{2}(x+1)$ by the maximality of $l'$,
which implies that $x\in 2A$ by Proposition \ref{pigeonhole}.
By symmetry, we have that $x+n\in 2A$ for any $0\leqslant x<n-r'$.
Hence $2A$ contains the interval $[l'+1,2n-r'-1]$. The length of the interval
is $2n-r'-l'-1$, which is greater than or equal to $2|A|-1$ because
\begin{eqnarray*}
\lefteqn{2n-r'-l'-1=2(n-r'-l'-1)+l'+1+r'+1-1}\\
& &\geqslant 2A(l'+1,n-r'-1)+2A(0,l')+2A(n-r',n)-1=2|A|-1.
\end{eqnarray*}

\medskip

\noindent {\bf Proof of Lemma \ref{lr}}\quad
Assume to the contrary that $l\geqslant n-r$. Clearly, $l\not=n-r$ by (\ref{bothsides}).
Hence we can assume that $l>n-r$. Let
\begin{equation}\label{definitionr}
r_0=\min\left\{x\in\left[n-l,r\,\right]:A(n-x,n)\leqslant\frac{1}{2}(x+1)
\right\}.
\end{equation}
By (\ref{bothsides}) we have that $n-r_0<l$. Let
\begin{equation}\label{definitionl}
l_0=\min\left\{x\in [n-r_0,l]:A(0,x)\leqslant\frac{1}{2}(x+1)\right\}.
\end{equation}
 We  have that $n-r_0<l_0$ again by (\ref{bothsides}). By the minimality of $l_0$ and $r_0$,
it is true that $A(0,x)>\frac{1}{2}(x+1)$ and $A(x,n)>\frac{1}{2}(n-x+1)$ for any $x\in
H\cap[n-r_0+1,l_0-1]$. So every hole in $[n-r_0+1,l_0-1]$ is an unstable hole.
Thus
\begin{equation}\label{middleh3}
H(n-r_0+1,l_0-1)\leqslant h_3.
\end{equation}
Now we have that
\begin{eqnarray}\label{allequal1}
\lefteqn{h\geqslant H(0,l_0)+H(n-r_0,n)-H(n-r_0,l_0)}\\
& &\geqslant\frac{1}{2}(l_0+1)+\frac{1}{2}(r_0+1)-H(n-r_0,l_0)\\
& &\geqslant\frac{1}{2}(n+1)+\frac{1}{2}(l_0-(n-r_0)+1)-H(n-r_0,l_0)\\
& &\geqslant\frac{1}{2}(|A|+h)-\frac{1}{2}H(n-r_0,l_0).
\end{eqnarray}
By solving the inequality above,
we get that $h\geqslant |A|-H(n-r_0,l_0)$, which implies that
\begin{equation}\label{allequal2}
0\geqslant |A|-2-h-H(n-r_0,l_0)+2\geqslant h_3-H(n-r_0+1,l_0-1)\geqslant 0
\end{equation}
by (\ref{wholeh3}) and (\ref{middleh3}).
Thus all inequalities in (\ref{allequal1})--(\ref{allequal2})
become equalities. In particular, it is true that
\begin{equation}\label{unstableholes}
H(n-r_0,l_0)=l_0-(n-r_0)+1=h_3+2.
\end{equation}
Notice that (\ref{unstableholes}) implies that
$[n-r_0,l_0]\cap A=\emptyset$ and the set of
all unstable holes is exactly the interval $[n-r_0+1,l_0-1]$.
Notice also that $l_0$ is a left stable hole and $n-r_0$ is a right stable hole.
These facts are important for the rest of the proof.

All arguments above this line are essentially due to Freiman in \cite{freiman2}.
The remaining part of the proof is new.

Notice that if $A(0,l_0)<\frac{1}{2}(l_0+1)$, then $2A(0,l_0-1)=2A(0,l_0)\leqslant l_0$,
which contradicts the minimality of $l_0$. Hence we have that
\begin{equation}\label{halfl}
A(0,l_0)=\frac{1}{2}(l_0+1).
\end{equation}
By the same reason, and the minimality of $r_0$, we have that
\begin{equation}\label{halfr}
A(n-r_0,n)=\frac{1}{2}(r_0+1).
\end{equation}
Since (\ref{halfl}), (\ref{halfr}), $l_0$ is left stable hole, and $n-r_0$ is a right
stable hole, we have, by Proposition \ref{addminim}, that $A\cap[0,l_0]$ and
and $A\cap[n-r_0,n]$ are anti-symmetric. Let
\begin{equation}\label{definitionab}
a=\max (A\cap[0,n-r_0])\,\mbox{ and }\,b=\min (A\cap[l_0,n]).
\end{equation}
Then $a<n-r_0$, $b>l_0$, and
$b-a\geqslant 2+l_0-(n-r_0)\geqslant 3$.
Since \[A(0,a)=A(0,l_0)=\frac{1}{2}(l_0+1)\geqslant\frac{1}{2}(a+3)\geqslant
\frac{1}{2}(a+1)+1,\] we have that $a>0$ and
\begin{equation}\label{gcd0a}
\gcd(A\cap[0,a])=1.
\end{equation}
By the same reason, we have that $b<n$ and
\begin{equation}\label{gcdbn}
\gcd(A\cap[b,n]-b)=1.
\end{equation}
By part 3 of Proposition \ref{sub3k-3prop}, we can assume that
$|2(A\cap[0,b])|\geqslant 3A(0,b)-3$. By the same reason, we can assume that
$|2(A\cap[a,n])|\geqslant 3A(a,n)-3$. By Proposition \ref{sub3k-3prop}, we have that
(\ref{intersection}) and (\ref{sub3k-3}) are true and the set in (\ref{twopieces})
is empty. We now use these facts to derive contradictions. Let
\[a'=\max A\cap[0,a-1]\,\mbox{ and }\,b'=\min A\cap[b+1,n].\]
A contradiction will be derived under each of the following conditions:

$b'-b<a-a'$,

$b'-b>a-a'$,

$b'-b=a-a'>b-a$,

$1<b'-b=a-a'\leqslant b-a$, and

$b'-b=a-a'=1$.

\medskip

Assume that $b'-b<a-a'$. Then we have the inequality
\[a'+b<a'+b'<a+b.\]

If $a'+b'=x+y$ such that $x,y\in A\cap[0,b]$ and $x\leqslant y$, then $y\not=b$ and
$y>a'$. This is true
because if $y=b$, then $x$ must be a number strictly between $a'$ and $b$ but not $a$,
and if $y\leqslant a'$ then $x+y\leqslant 2a'<a'+b'$.
Hence the only possible choice for $y$ is $a$. With $y=a$, we have that $x=a$ because
$x\leqslant a'$ implies $x+y\leqslant a'+a<a'+b'$.  So $a'+b'=2a$.

If $a'+b'=x+y$ such that $x,y\in A\cap[a,n]$ and $x\leqslant y$, then again we have that
$a'+b'=2a$ because $a'+b'<a+b$.

Since the set in (\ref{twopieces}) is empty, we have that
\[a'+b'\in (2(A\cap[0,b]))\cup (2(A\cap[a,n])),\]
which implies that $a'+b'=2a$ by the arguments above. As a consequence,
we have that $2a\not\in b+A\cap [0,b]$.

The fact that $2a\not\in b+A\cap [0,b]$
 will be used in the next several paragraphs
to show that $|2(A\cap [0,b])|>3A(0,b)-3$, which contradicts (\ref{sub3k-3}).

Let $z=\max\{x\in [-1,l_0-1]:A(0,x)\leqslant\frac{1}{2}(x+1)\}$. Clearly, $z+1\in A$
by the maximality of $z$. Notice that
$A(0,z)=\frac{1}{2}(z+1)$ and $A\cap [z+1,l_0]$ is left dense in $[z+1,l_0]$
by the maximality of $z$ and (\ref{halfl}).

If $z=-1$, then $A\cap[0,l_0]$ is left dense in $[0,l_0]$ and hence
\[2(A\cap [0,b])\supseteq [0,l_0-1]\cup(b+A\cap [0,b])\cup\{2a\}.\]
So $|2(A\cap [0,b])|\geqslant 3A(0,b)-2>3A(0,b)-3$,
which contradicts (\ref{sub3k-3}).

Suppose that $z>-1$. Then $z>0$ because $A(0,0)=1>\frac{1}{2}$.

If $\gcd(A\cap[0,z+1])=1$,
then $|2(A\cap [0,z+1])|\geqslant 3A(0,z+1)-3$
by part 3 of Proposition \ref{sub3k-3prop}.
If $z\in A$, then $z+b\in 2A$ and
\begin{eqnarray*}
\lefteqn{|2(A\cap [0,b])|}\\
& &\geqslant |2(A\cap[0,z+1])|-1\\
& &\quad +|[2z+2,z+l_0]|+|(b+A\cap [z,b])\cup\{2a\}|\\
& &\geqslant
3A(0,z+1)-4+l_0-z-1+A(z,b)+1\\
& &\geqslant 3A(0,z)-1+2A(z+1,l_0)+A(z+1,b)+1\\
& &=3A(0,z)+3A(z+1,b)-2>3A(0,b)-3,
\end{eqnarray*} which contradicts (\ref{sub3k-3}).
So we can assume that $z\not\in A$. Let \[z'=\max A\cap [0,z-1].\]
Then $z'+z+2\in (2A)\smallsetminus (2(A\cap[0,z+1]))$.
 Hence
 \begin{eqnarray*}
 \lefteqn{|2(A\cap [0,b])|}\\
 & &\geqslant |(2(A\cap[0,z+1]))\cup\{z'+z+2\}|-1\\
& &\quad +|[2z+2,z+l_0]|+|(b+A\cap [z+1,b])\cup\{2a\}|\\
& &\geqslant
3A(0,z+1)-3+l_0-z-1+A(z+1,b)+1\\
& &\geqslant 3A(0,z)+2A(z+1,l_0)+A(z+1,b)\\
& &=3A(0,z)+3A(z+1,b)-2>3A(0,b)-3,
\end{eqnarray*} which again contradicts (\ref{sub3k-3}).

Thus we can assume that $\gcd(A\cap [0,z+1])=d>1$. Clearly, $d=2$ and
$A\cap [0,z+1]$ is an arithmetic progression of
difference $2$ by the fact that $A(0,z)=\frac{1}{2}(z+1)$. Hence
\begin{eqnarray*}
\lefteqn{|2(A\cap[0,b])|}\\
& &\geqslant |A\cap[0,z-1]+A\cap[0,z+1]|+|(z+2)+A\cap[0,z-1]|\\
& &\quad +|[2z+2,z+l_0]|+|(b+A[z+1,b])\cup\{2a\}|\\
& &\geqslant 3A(0,z-1)+l_0-z-1+A(z+1,b)+1\\
& &\geqslant 3A(0,z)+2A(z+1,l_0)+A(z+1,b)\\
& &=3A(0,b)-2>3A(0,b)-3,
\end{eqnarray*} which again contradicts
(\ref{sub3k-3}).

Assume that $a-a'<b'-b$. The proof is symmetric to the case for $a-a'>b'-b$.

Assume that $b'-b=a-a'=d'$.

If $d'>b-a$, then
$2a\not\in b+A\cap[0,b]$. Hence $|2(A\cap [0,b])|>3A(0,b)-3$
by the same argument as above, which contradicts (\ref{sub3k-3}).
Thus, we can now assume that $d'\leqslant b-a$.

Suppose that $1<d'\leqslant b-a$. Let $a''$ be the greatest element
in $A\cap[0,a']$ which is not congruent to $a$ modulo $d'$. The number $a''$ exists by
(\ref{gcd0a}).

If $b'+a''\in 2(A\cap[a,n])$, then $b'+a''=2a$ because
$b'+a''<a+b$. This implies
that $2a=b+(d'+a'')\not\in b+A\cap[0,b]$ by the maximality of $a''$. Hence
$|2(A\cap [0,b])|>3A(0,b)-3$, which contradicts (\ref{sub3k-3}).
Since the set in (\ref{twopieces}) is empty,
we can assume that $b'+a''\in 2(A\cap[0,b])$.

Clearly, $b'+a''\not\in b+A\cap[0,b]$ by the maximality of $a''$.
Let $b'+a''=x+y$ for some $x,y\in A\cap[a''+1,a]$. Then $b'+a''\leqslant 2a$
and $x,y$ are congruent to $a$ modulo $d'$.
Since the fact that
$2a\not\in b+A\cap [0,b]$ contradicts (\ref{sub3k-3}), we can assume that
$2a=b+z$ for some $z\in A\cap[0,a]$. Hence $d'+a''=b'-b+a''\leqslant 2a-b=z$.
By the maximality of $a''$, we have that $z$ is congruent to $a$ modulo $d'$.
This implies that that $b$ is congruent to $a$ modulo $d'$
because $2a=b+z$. Hence $b'$ is congruent to $a$ modulo $d'$.
Now we have that $a''=x+y-b'$ is congruent to
$a$ modulo $d'$, which contradicts the definition of $a''$.

We can now assume that $d'=1$, i.e.,
\begin{equation}\label{a-1andb-1}
a'=a-1\in A\,\mbox{ and }\,b'=b+1\in A.
\end{equation}
The derivation of a contradiction under this case is much harder that the previous
cases. Notice that $A$ also satisfies the condition (\ref{a-1andb-1})
when $A$ is a bi-arithmetic progression of difference $4$ such as $A=\{0,1,4,5,8,9\}$.

Since $A\cap[0,l_0]$ and $A\cap[n-r_0,n]$ are anti-symmetric
and $[n-r_0,l_0]\cap A=\emptyset$,
we have that $[0,l_0-(n-r_0)]\subseteq A$ and $[2n-l_0-r_0,n]\subseteq A$.
In particular, we have that
\begin{equation}\label{01n-1n}
0,1,n-1,n\in A.
\end{equation}

Next we prove four claims for the existence of unstable holes if $A$ has a certain configuration. These claims will be used to derive a contradiction.

\begin{claim}\label{claim1}
If $z\in A$, then $z-1\in A$ or $z+1\in A$.
\end{claim}

\noindent {\bf Proof of Claim \ref{claim1}}\quad Suppose that $z-1,z+1\not\in A$ and
$z\in A$. We call such an element $z$ an isolated point of $A$. Then
$z\in [3,n-3]$ and $z\not\in [a-2,b+2]$ by (\ref{a-1andb-1}) and (\ref{01n-1n}).
If $A(0,z-1)>\frac{1}{2}z$, then $z-1\in 2A$ and
$n+z-1=(n-1)+z\in 2A$. Hence $z-1$ is an unstable hole,
which contradicts (\ref{unstableholes}). Hence $A(0,z-1)\leqslant\frac{1}{2}z$,
which implies that $A(0,z+1)=A(0,z-1)+1\leqslant\frac{1}{2}z+1=\frac{1}{2}(z+2)$.
By (\ref{bothsides}) we have that $A(z+1,n)>\frac{1}{2}(n-z)$. Therefore,
$n+z+1\in 2A$ by Proposition \ref{pigeonhole} and
$z+1\in 2A$ because $z\in A$ and $1\in A$, which again
contradicts (\ref{unstableholes}).

\medskip

Claim \ref{claim1} says that $A$ does not contains any isolated points in $A$.

\begin{claim}\label{claim2}
If $z\in H$, then either $z-1\in H$ or $z+1\in H$.
\end{claim}

\noindent {\bf Proof of Claim \ref{claim2}}\quad Suppose that
$z-1,z+1\in A$ and $z\in H$. Since $l_0,n-r_0\not\in A$, we have that
$z\not\in [n-r_0+1,l_0-1]$. Since $z=z-1+1\in 2A$ and $z+n=z+1+n-1\in 2A$,
it is true that $z$ is an unstable hole, which contradicts (\ref{unstableholes}).

\medskip

Claim \ref{claim2} says that there do not exist any isolated holes of $A$.

\begin{claim}\label{claim3}
(a) If $0<x<y<z<n$ are such that $x,z,z+1\in H$, $y\in A$, and $A(0,z)=\frac{1}{2}(z+1)$,
then $z+1$ is an unstable hole.

(b) If $0<x<y<z<n$ are such that $x-1,x,z\in H$, $y\in A$, and $A(x,n)=\frac{1}{2}(n-x+1)$,
then $x-1$ is an unstable hole.
\end{claim}

\noindent {\bf Proof of Claim \ref{claim3}}\quad We prove (a) only
and (b) follows by symmetry.
Without loss of generality, let
$x=\max H\cap[0,y]$. By (\ref{bothsides}) we have that $A(z+1,n)=A(z,n)>\frac{1}{2}(n-z+1)
>\frac{1}{2}(n-z)$. Hence $n+z+1\in 2A$. Notice that $A(0,z)=\frac{1}{2}(z+1)$ implies that
$z\not\in [n-r_0+1,l_0-1]$. So $z$ is not an unstable hole by (\ref{unstableholes}).
Since $z+n\in 2A$ by Proposition \ref{pigeonhole}, we have that
$z\not\in 2A$. By part 1 of Proposition \ref{addminim},
$A\cap [0,z]$ is anti-symmetric in $[0,z]$. So $x\in z-A\cap[0,z]$.
Hence $x+1\in A\cap[0,z]\cap(z+1-A\cap[0,z])$, which implies that
$z+1\in 2A$ and hence $z+1$ is an unstable hole.

\medskip

Claim \ref{claim3} (a) implies that if $[0,a]\not\subseteq A$, then $b=l_0+1$ because $b>l_0+1$ implies that
$l_0+1$ is an unstable hole, which contradicts (\ref{unstableholes}).
By symmetry, Claim \ref{claim3} (b) implies that $a=n-r_0-1$ if $[b,n]\not\subseteq A$.

\begin{claim}\label{claim4}
If $[x,y]\subseteq H$ is a gap of $A$ with $y-x\geqslant 2$,
$H\cap[0,x-1]\not=\emptyset$, and $H\cap[y+1,n]\not=\emptyset$,
then $[x,y]$ contains an unstable hole.
\end{claim}

\noindent {\bf Proof of Claim \ref{claim4}}\quad If $A(0,x)\leqslant\frac{1}{2}(x+1)$, then
$A(x,n)>\frac{1}{2}(n-x+1)$ by (\ref{bothsides}), which implies that $n+x\in 2A$
by Proposition \ref{pigeonhole}.
Also $x=x-1+1\in 2A$. Hence $x$ is an unstable hole.
Symmetrically, if $A(y,n)\leqslant\frac{1}{2}(n-y+1)$,
then $y\in 2A$ and $n+y=(n-1)+(y+1)\in 2A$. Hence $y$ is an unstable hole.
So we can now assume that $A(0,x)>\frac{1}{2}(x+1)$ and $A(y,n)>\frac{1}{2}(n-y+1)$.

Let $t\in [x+1,y-1]$, $x'=\max H\cap[0,x-1]$, and $y'=\min H\cap[y+1,n]$.
If $A(0,t)>\frac{1}{2}(t+1)$ and $A(t,n)>\frac{1}{2}(n-t+1)$, then
$t$ is an unstable hole. Otherwise we can assume, without loss of generality, that
$A(0,t)\leqslant\frac{1}{2}(t+1)$.

Assume, without loss of generality again, that
$t\in [x+1,y-1]$ is the least element such that
$A(0,t)\leqslant\frac{1}{2}(t+1)$. Notice that
$t+1\leqslant y$ and $t+1\in H$. If $A(0,t)<\frac{1}{2}(t+1)$, then
$2A(0,t-1)=2A(0,t)\leqslant t$. Hence $A(0,t-1)\leqslant\frac{1}{2}t$.
Since $A(0,x)>\frac{1}{2}x$, we have that $x<t-1$. This contradicts the
minimality of $t$. Therefore, we can assume that $A(0,t)=\frac{1}{2}(t+1)$.
Now we conclude that $t+1$ is an unstable hole by Claim \ref{claim3}.

\medskip

Claim \ref{claim4} says that if $A$ has a gap $[x,y]$ of length at least $3$, i.e.,
$y-x\geqslant 2$,
then $[x,y]$ is either the first gap or the last gap or the middle gap $[a+1,b-1]$ of $A$.

\medskip

We now continue the proof of Theorem \ref{main1}
by deriving a contradiction under the assumption that $d'=1$, i.e.,
$a-1,a,b,b+1\in A$.

If $n-b<b-a$ and $a<b-a$, then $A$ is a subset of the bi-arithmetic progression
$[0,a]\cup [b,n]$ of difference $1$. So $|2A|=3|A|-3$ implies that $A=[0,a]\cup[b,n]$
by Theorem \ref{2k-1+b}. Hence part 1 of Theorem \ref{main1} is true.
Thus we can now assume that either $n-b\geqslant b-a$ or $a\geqslant b-a$.

Without loss of generality let $a\geqslant b-a$. If $A\cap[0,a]=[0,a]$, then
$A(0,l_0)=A(0,a)=a+1\geqslant\frac{1}{2}b+1\geqslant\frac{1}{2}(l_0+1)+1=A(0,l_0)+1$,
 which is absurd. So we can assume that $H\cap[0,a]\not=\emptyset$.
Let
\begin{equation}\label{definitionz}
z=\min\left\{x\in [0,l_0]:A(0,x)\leqslant\frac{1}{2}(x+1)\right\}.
\end{equation}
Since $A(0,l_0)=\frac{1}{2}(l_0+1)$, the number $z$ is well defined and $z\not=a$.
We now divide the rest of the proof into two cases: $z>a$ or $z<a$. In each case
we derive a contradiction.

\medskip

\noindent {\bf Case 1}\quad $z>a$.

\medskip

We want to show that $|2(A\cap[0,b])|\geqslant 3A(0,b)-2$,
which contradicts (\ref{sub3k-3}).

Notice that $z=l_0$ because $[a+1,l_0]\cap A=\emptyset$. So $A\cap [0,l_0]$ is
left dense, anti-symmetric in $[0,l_0]$.
Let $y=\min H\cap[0,a]$. If $y+b\not=2a$,
then $y+b=(y-1)+(b+1)\in (2A)\smallsetminus (2(A\cap[a,n]))$,
which implies that $y+b\in 2(A\cap[0,b])$ by the fact that the set in (\ref{twopieces})
is empty. If $y+b=2a$, then again $y+b\in 2(A\cap[0,b])$.
Hence
\begin{eqnarray*}
\lefteqn{|2(A\cap[0,b])|}\\
& &\geqslant |[0,l_0-1]|+|(b+A\cap[0,b])\cup\{y+b\}|\\
& &\geqslant 2A(0,l_0)-1+A(0,b)+1=3A(0,b)-2.
\end{eqnarray*}

\medskip

\noindent {\bf Case 2}\quad $z<a$.

\medskip

The proof of this case is much longer than the proof of Case 1.

Notice that $z\not\in A$ by the minimality of $z$ and $z>2$.
If $z-1\in A$, then $A(0,z-2)=A(0,z)-1\leqslant
\frac{1}{2}(z+1)-1=\frac{1}{2}(z-1)$, which contradicts
the minimality of $z$. Hence $z-1\not\in A$. Notice that $A(0,z-1)>\frac{1}{2}z$
by the minimality of $z$. If $A(z-1,n)>\frac{1}{2}(n-z+2)$, then $z-1$ is an unstable
hole below $a$ by Proposition \ref{pigeonhole},
which contradicts (\ref{unstableholes}). Hence we can assume that
$A(z-1,n)\leqslant\frac{1}{2}(n-z+2)$. Since
\begin{eqnarray*}
\lefteqn{A(z-1,n)=|A|-A(0,z-2)}\\
& &=|A|-A(0,z)\geqslant\frac{1}{2}(n+3)-\frac{1}{2}(z+1)
=\frac{1}{2}(n-z+2),
\end{eqnarray*}
we have that
\begin{equation}\label{rightside}
A(z-1,n)=\frac{1}{2}(n-z+2).
\end{equation}
By Claim \ref{claim3} (b), we can assume that $z-2\in A$ because
otherwise $z-2$ becomes an unstable hole below $a$.
It is worth mentioning
that (\ref{rightside}) and $A(0,z)\leqslant\frac{1}{2}(z+1)$ imply
\begin{equation}\label{hplus2}
|A|=\frac{1}{2}(n+3)=\frac{1}{2}(|A|+h+2),
\end{equation}
which implies $|A|-2=h$ and
\begin{equation}\label{nounstablehole}
h_3=|A|-2-h=0.
\end{equation}
So $A$ has no unstable holes and $n-r_0=l_0-1$.
We now divide the rest of the proof
into two cases: $z>3$ or $z\leqslant 3$.

\medskip

{\bf Case 2.1}\quad $z>3$.

\medskip

If $z-3\not\in A$, then $z-2\in A$ is an isolated point in $A$,
which contradicts Claim \ref{claim1}. But if $z-3\in A$, then
\[A(0,z-4)=A(0,z)-A(z-3,z)\leqslant\frac{1}{2}(z+1)-2=\frac{1}{2}(z-3),\]
which contradicts the minimality of $z$.

\medskip

{\bf Case 2.2}\quad $z\leqslant 3$.

\medskip

Since $A(0,z)=\frac{1}{2}(z+1)$, $z\not\in A$, and $0,1\in A$, we have that $z=3$ and
$z-1=2\not\in A$. Hence $A\cap[0,3]=\{0,1\}$. Notice that we have assumed that
$A$ is not a bi-arithmetic progression (of difference $4$).

Let \[V=\{x\in [0,n]:x\equiv 0,1\,(\mod 4)\}.\]
Then $A\not=V$.
Let \[z'=\min\{x\in [0,n]:A\cap[0,x]\not=V\cap[0,x]\}.\]
 Notice that $n\geqslant z'>z=3$ and $A\cap[0,z'-1]=V\cap[0,z'-1]$ is the maximal bi-arithmetic
progression of difference $4$ inside $A$ containing $0,1$.
We now divide the rest of the proof
into four cases in terms of the value of $z'$ modulo $4$.

\medskip

\quad {\bf Case 2.2.1}\quad $z'\equiv 0\,(\mod 4)$.

\medskip

Clearly, $z'\not\in A$ because otherwise $A\cap[0,z']=V\cap[0,z']$.

If $z'>4$, then $A\cap[0,z']=\{0,1,4,5,\ldots,z'-4,z'-3\}$ and $z'$ is at least $8$.
Since $A(0,z'-1)=\frac{1}{2}z'$ by the definition
of $V$, $3,z'-1,z'\in H$, and $4\in A$, we have that $z'$ is an unstable hole
by Claim \ref{claim3}, which contradicts (\ref{nounstablehole}).

So we can now assume that $z'=4$, which implies that $A\cap[0,4]=\{0,1\}$.

Let $c=\min A\cap[z',a]$.

Recall that $l_0$ is a left stable hole and
$A\cap [0,l_0]$ is anti-symmetric in $[0,l_0]$ by Proposition \ref{addminim}.
Since $0,1,c\in A$ and $2\not\in A$,
we have that $l_0,l_0-1,l_0-c\not\in A$ and $l_0-2\in A$.
Consequently, $l_0-2=a$ by (\ref{definitionab}) and $l_0-1=n-r_0$.

Suppose that $H\cap[c,a]\not=\emptyset$ and
let $t=\max H\cap[c,a]$. Then $t+1\in A$ and
$t-1\not\in A$ by Claim \ref{claim2}. Since again
$A\cap[0,l_0]$ is anti-symmetric in $[0,l_0]$, we have that
$l_0-[2,c-1]=[l_0-c+1,a]\subseteq A$ because $[2,c-1]\subseteq H$.
Consequently, $t=l_0-c$.

Since $c-2>2$, we have that $a-t>2$, which implies that
$A(t,a)=A(t+1,a)=a-t>\frac{1}{2}(a-t+2)$. Hence
\begin{eqnarray*}
\lefteqn{A(t-1,n)=A(t,n)=A(t,a)+A(n-r_0,n)}\\
& &>\frac{1}{2}(a-t+2)+\frac{1}{2}(r_0+1)=\frac{1}{2}(a-t+2)+\frac{1}{2}(n-l_0+2)\\
& &=\frac{1}{2}(a-t+2)+\frac{1}{2}(n-a)=\frac{1}{2}(n-t+2)>\frac{1}{2}(n-t+1),
\end{eqnarray*}
which implies that $n+t-1\in 2A$. If $t-2\not\in A$, then the gap containing $t$ has
length $\geqslant 3$, which implies that the gap contains an unstable hole
by Claim 4. Hence we have a contradiction to (\ref{nounstablehole}).
Therefore, we can now assume that $t-2\in A$. But this implies that $t-1=(t-2)+1\in 2A$. So
$t-1$ is unstable hole, which again contradicts (\ref{nounstablehole}).

We can now assume that $H\cap[c,a]=\emptyset$, i.e., $A\cap[0,b]=\{0,1\}\cup [c,a]\cup\{b\}$.
We want to show that $|2(A\cap[0,b])|\geqslant 3A(0,b)-2$, which contradicts
(\ref{sub3k-3}).

Since $A\cap[0,l_0]$ is anti-symmetric in $[0,l_0]$, we have that
$a-c=c-3\geqslant 2$. Notice that $b=l_0+1=a+3$ and
$A(0,b)=a-c+4$. Now we have that
\[2(A\cap[0,b])\supseteq [0,2]\cup[c,a+1]\cup[2c,a+b]
\cup\{2b\}\,\mbox{ and}\]

\vspace{-0.4in}

\begin{eqnarray*}
\lefteqn{|2(A\cap[0,b])|\geqslant 3+a-c+2+a+b-2c+1+1}\\
& &=3a-3c+10=3A(0,b)-12+10=3A(0,b)-2.
\end{eqnarray*}

\medskip

\quad {\bf Case 2.2.2}\quad $z'\equiv 1\,(\mod 4)$.

\medskip

We have that $z'\not\in A$, $z'-1\in A$, and $z'-2\not\in A$.
Hence $z'-1$ is an isolated point of $A$, which contradicts Claim \ref{claim1}.

\medskip

\quad {\bf Case 2.2.3}\quad $z'\equiv 2\,(\mod 4)$.

\medskip

We have that $z',z'-1,z'-2\in A$ and $z'-3\not\in A$.

Let $c=\max\{x\geqslant z':[z',x]\subseteq A\}$.

Notice that $[z'-2,c]\subseteq A$. We divide the proof of this case into
four subcases for $c=n$, $b<c<n$, $c=a$, or $c<a$. Notice that $c=b$ is impossible
because $c-1\in A$.

\medskip

\quad\quad {\bf Case 2.2.3.1}\quad $c=n$.

\medskip

Since $A=(V\cap[0,z'-3])\cup [z'-2,n]$, we have that
\begin{eqnarray*}
\lefteqn{|A|=A(0,z'-3)+A(z'-2,n)}\\
& &=\frac{1}{2}(z'-2)+n-z'+3\\
& &=\frac{1}{2}(n+1)+\frac{1}{2}(n-z'+3)\\
& &\geqslant\frac{1}{2}(n+1)+\frac{3}{2}=\frac{1}{2}(n+4),
\end{eqnarray*}
which contradicts (\ref{hplus2}).

\medskip

\quad\quad {\bf Case 2.2.3.2}\quad $b<c<n$.

\medskip

Recall that $a+3=b$. Since $z',z'-1,z'-2\in A$, we have that $b\leqslant z'-2$.
Let $x=2n-r_0-c$ and $y=2n-r_0-z'+2$. Notice that $z'-3\not\in A$, $[z'-2,c]\subseteq A$,
and $c+1\not\in A$.
Since $n-r_0$ is a right stable hole and $A\cap [n-r_0,n]$ is anti-symmetric in
$[n-r_0,n]$, we have that $[x,y]$ is a gap of $A$ with length
$y-x+1=c-z'+3\geqslant 3$.
Notice that $c<x$ because gaps of $A$ below
$c$ are also gaps of $V$ with length $2$ while the length of $[x,y]$ is at least $3$.

Suppose that $c+1<x$. Since $2n-r_0-[c+1,x-1]=[c+1,x-1]$ and $A\cap [n-r_0,n]$
is anti-symmetric in $[n-r_0,n]$, we
 have that $A(c+1,x-1)=\frac{1}{2}(x-c-1)$.
Let \[c'=\max\{t\in [c+1,x-1]:[c+1,t]\subseteq H\}.\] Since $x-1\in A$,
we have that $c'<x-1$.
If $A(0,c')>\frac{1}{2}(c'+1)$, then $c'\in 2A$ and $n+c'=(n-1)+(c'+1)\in 2A$,
which contradicts (\ref{nounstablehole}). Hence we can assume that $A(0,c')\leqslant
\frac{1}{2}(c'+1)$. Since $A(0,z'-3)=\frac{1}{2}(z'-2)$, we have that
$A(z'-2,c')=A(z'-2,c)=c-z'+3\leqslant\frac{1}{2}(c'-z'+3)$, which implies that $c'-c\geqslant
c-z'+3\geqslant 3$. Hence $[c+1,c']$ is a gap of $A$ with length at least $3$.
Since $l_0<b<c$ and $c'<x-1<x$, the gap $[c+1,c']$ contains an unstable hole by Claim \ref{claim4}, which again contradicts (\ref{nounstablehole}).

Thus we can assume that $c+1=x$.

If $H\cap[y+1,n]\not=\emptyset$, then $[x,y]$ contains an unstable hole by Claim \ref{claim4},
which contradicts (\ref{nounstablehole}).
So we can assume that $H\cap[y+1,n]=\emptyset$,
which means that $y+1=n-1$ because $2n-r_0-(n-2)=n-r_0+2=b\in A$
implies $n-2\not\in A$. Hence
\[A=(V\cap[0,z'-3])\cup [z'-2,c]\cup\{n-1,n\}.\]
Since $A\cap[n-r_0,n]$ is anti-symmetry in $[n-r_0,n]$, we have that $(n-2)-(c+1)=c-(z'-2)$
or $n-c-3=c-z'+2$. Notice also that $a+3=b=z'-2$.
We are now ready to show that $|2(A\cap[a,n])|\geqslant 3A(a,n)-2$, which will
contradicts (\ref{sub3k-3}).

Notice that $A\cap[a,n]=\{a\}\cup [b,c]\cup\{n-1,n\}$ and $A(a,n)=c-b+4$.
Since \[2(A\cap[a,n])=\{2a\}\cup[a+b,2c]\cup[n-1+b,n+c]\cup[2n-2,2n],\] we have that
\begin{eqnarray*}
\lefteqn{|2(A\cap[a,b])|=1+2c-a-b+1+c-b+2+3}\\
& &=3c-3b+10=3A(a,n)-12+10=3A(a,n)-2.
\end{eqnarray*}

\quad\quad {\bf Case 2.2.3.3}\quad $c=a$.

\medskip

Since $A\cap[0,l_0]$ is anti-symmetric in $[0,l_0]$, we have that
$[x,y]=l_0-[z'-2,a]=[2,l_0-z'+2]\subseteq H$ is a gap of $A$ with length
at least $3$, which implies that $4\not\in A$.
But from the first paragraph of Case 2.2, we have that $z'\geqslant 4$ and from
the assumption of Case 2.2.3, i.e., $z'\equiv 2\,(\mod 4)$, we have that
$z'\geqslant 6$. Thus we have that $4\in A$ by the definition of $z'$.
So we have a contradiction.

\medskip

\quad\quad {\bf Case 2.2.3.4}\quad $c<a$.

\medskip

Notice again that $z'\geqslant 6$ because $z'>3$ and $z'\equiv 2\,(\mod 4)$.
Hence $0,1,4,5\in A$ and $2,3\not\in A$. Since $A\cap[0,l_0]$ is anti-symmetric
in $[0,l_0]$, we have that $a,a-1\in A$ and $a-2,a-3\not\in A$.
Since $[z'-2,c]\subseteq A$ and $z'-3,c+1\not\in A$, there is
a gap $[x,y]=[l_0-c,l_0-z'+2]\subseteq H$ of $A$ with length at least $3$.
Notice that $[x,y]$ is not the gap $[2,3]$ of length $2$. We have that
$x\geqslant 6$. Clearly, $y<a-1$. By Claim \ref{claim4},
$[x,y]$ contains an unstable hole, which contradicts (\ref{nounstablehole}).

\medskip

\quad {\bf Case 2.2.4}\quad $z'\equiv 3\,(\mod 4)$.

\medskip

By the definition of $z'$ we have that $z',z'-2\in A$ and $z'-1\not\in A$.
Therefore, $z'-1$ is an isolated hole, which contradicts Claim \ref{claim2}.

This completes the proof of Theorem \ref{main1} when $|A|>\frac{1}{2}(n+2)$.

\begin{remark}
Theorem \ref{2k-1+b} and Theorem \ref{3k-3} characterize the structure of $A$ when
$|2A|\leqslant 3|A|-3$. The structure of $A$ in Theorem \ref{freimanthm}
 is given indirectly by describing a property for $2A$, e.g. $2A$ contains
 an interval of length $2|A|-1$.
 In fact, Freiman's original result in \cite{freiman2}
shows the following: If $e$ is the greatest $x$ in $[-1,n]$ such that
$x\not\in 2A$ and $c$ is the least $x$ in $[0,n+1]$ such that $x+n\not\in 2A$, then
$e<c$, which implies that $2A$ contains an interval of length $2|A|-1$.
Hence the structural information in \cite{freiman2} is presented directly for $A$
instead of $2A$.

Assume $|2A|<3|A|-3$. Let $l'$ and $r'$ be the maximal $l$ and $r$, respectively, as
defined in Lemma \ref{lr}. Then $e\leqslant l'$ and $c\geqslant n-r'$ by
Proposition \ref{pigeonhole} and the maximality of $l'$ and $r'$.
It could happen that $e$ is strictly less than $l'$ and
$c$ is strictly greater than $n-r'$.
It is not too difficult to modify Freiman's proof to
show that $l'<n-r'$, which implies that $e<c$.
As an extra benefit, the conclusion $l'<n-r'$
gives some geometric information directly for $A$.
Roughly speaking, $l'<n-r'$
indicates that $A$ is thin in $[0,l']$ and in $[n-r',n]$, and $A$ is thick
in $[l'+1,n-r'-1]$.

By the comments above, we can say that Lemma \ref{lr} is slightly better than
the statement that $\,2A$ contains an interval of length $2|A|-1$
when $|A|\geqslant\frac{1}{2}(n+3)$.
\end{remark}

\subsection{Proof of Theorem \ref{main1} when $|A|=\frac{1}{2}(n+2)$}

Throughout this subsection we assume that
\begin{equation}\label{least}
|A|=\frac{1}{2}(n+2).
\end{equation}

Notice that (\ref{least}) cannot occur when $n$ is an odd number.

Let $x$ be a hole in $A$. We call $x$ a {\it balanced hole} if $A(0,x)=\frac{1}{2}(x+1)$
and $A(x,n)=\frac{1}{2}(n-x+1)$. Notice that if $A(0,y)=\frac{1}{2}(y+1)$
and $A(y,n)=\frac{1}{2}(n-y+1)$ for some $y\in [0,n]$, then $y\not\in A$
and if $A(0,x)=\frac{1}{2}(x+1)$ for some hole $x$ in $A$, then
$x$ is a balanced hole by (\ref{least}).

We want to show that $A$ is in the form of either $T_{k,n}$ or $n-T_{k,n}$
for some $k\in [0,\frac{1}{2}n-2]$ where $T_{k,n}$ is defined in
(\ref{tk}) or is in the form of $S_{u,n}$ for some $u\in[4,n-6]$
where $S_{u,n}$ is defined in (\ref{su}). It is worth mentioning that if $n=10$
and $B$ is a Freiman isomorphism image of $K_6$ in (\ref{sixpoints}),
then $|B|=\frac{1}{2}(n+2)$ and $B=B_i$ for $i=1,2,3$, or $4$ where $B_i$'s are
defined in part 4 of Proposition \ref{sixpoints2}.
Notice that $B_1=T_{0,10}$, $B_2=T_{2,10}$, $B_3=n-T_{0,10}$, and $B_4=n-T_{2,10}$.

\medskip

\noindent {\bf Case 1}\quad $0,1\in A$.

\medskip

In this case we want to show that $A$ is in the form of $T_{0,n}$ or
$A$ is an arithmetic progression of difference $1$ or $4$.
Since we have assumed that $A$ is not a bi-arithmetic progression,
the latter is a contradiction. Let
\begin{equation}\label{balancedz}
z=\min\left\{x\in [0,n]:A(0,x)\leqslant\frac{1}{2}(x+1)\right\}.
\end{equation}

Since $0,1\in A$, we have that $z\geqslant 3$. Clearly, $z-1,z\not\in A$, because otherwise
$A(0,z-2)\leqslant\frac{1}{2}z$, which contradicts the minimality of $z$.
We also have that $A(0,z)=\frac{1}{2}(z+1)$ and $A\cap[0,z]$ is left dense
in $[0,z]$ by the minimality of $z$. Notice that $z$ is a balanced hole.

If $z=2|A|-3$, then $2A\supseteq [0,z-1]\cup(n+A)$. Hence \[3|A|-3=|2A|
\geqslant |[0,z-1]\cup(n+A)|=2|A|-3+|A|=3|A|-3,\] which implies that
$2A=[0,z-1]\cup(n+A)$. So $A$ is in the form of $T_{0,n}$.
Therefore, we can now assume that $z<2|A|-3=n-1$.

We now intend to derive a contradiction by showing that
either $|2A|>3|A|-3$ or $A$ is a bi-arithmetic progression
of difference $1$ or $4$.

Let $a=\max A\cap[0,z]$ and $b=\min A\cap[z,n]$. By part 3 of
Proposition \ref{sub3k-3prop}, we can assume that
\begin{equation}\label{sub3k-3-3}
|2(A\cap[0,b])|\geqslant 3A(0,b)-3.
\end{equation}
Since $A(z,n)=A(b,n)=\frac{1}{2}(n-z+1)\geqslant\frac{1}{2}(z+2-z+1)>1$, the set
$A\cap[a,n]$ contains at least three elements.

Suppose that $\gcd(A\cap [b,n]-b)>1$. Then $A(z,n)=\frac{1}{2}(n-z+1)$ implies that
$b=z+1$ and $A\cap [b,n]$ is an arithmetic progression of difference $2$.
Notice that
\begin{equation}\label{atleast}
2A\supseteq [0,z-1]\cup(A\cap[b,n-2]+\{0,1\})\cup (n+A).
\end{equation}
\[\mbox{So }\,|2A|\geqslant 2A(0,a)-1+2A(b,n-2)+|A|=3|A|-3.\]
Let $A\cap [0,a]=E_0\cup O_0$ where $E_0$ contains all even numbers and $O_0$
contains all odd numbers in $A\cap [0,a]$. If $E_0$ is not a set of consecutive
even numbers, let $x>0$ be such that $x\not\in E_0$ and $x+2\in E_0$.
Then $n+x=(n-2)+(x+2)$ is in $2A$ but not in the right side of (\ref{atleast}).
So we have $|2A|>3|A|-3$.
By the same reason we can assume that $O_0$ contains consecutive odd numbers.
If $a=z-1=b-2$, then $a$ is an even number because $n$ is even and $A\cap [b,n]$
is an arithmetic progression of difference $2$.
Since $A(0,z)=\frac{1}{2}(z+1)=|E_0|$, then
$O_0=\emptyset$, which contradicts $1\in A$. Hence we can assume that
$a<b-2$. So $b-2\not\in A$. Since $b+(n-2)\in 2A$ and $b-2+n
\not\in n+A$, we have that $|2A|>3|A|-3$ by (\ref{atleast}).
Notice that in the proof of $|2A|>3|A|-3$ above
when $b<n$ and $\gcd(A\cap[b,n]-b)>1$, we have never
tried to force $z\in 2A$. So if we can show that $z$ is in $2A$
when $\gcd(A\cap[b,n]-b)>1$, then
$|2(A\cap[0,b])|>3A(0,b)-3$. This fact will be used later.
We can now assume that $\gcd(A\cap[b,n]-b)=1$.

By part $3\,$ of Proposition \ref{sub3k-3prop},
we can assume that $|2(A\cap[a,n])|\geqslant 3A(a,n)-3$.
Together with (\ref{sub3k-3-3}), we can derive the same equalities as in
(\ref{intersection}), (\ref{sub3k-3}), and the set in (\ref{twopieces}) is empty
by Proposition \ref{sub3k-3prop}. As a consequence, we have that $A\cap[0,z]$
is left dense and anti-symmetric in $[0,z]$.

\medskip

{\bf Case 1.1}\quad $H\cap[0,a]=\emptyset$.

\medskip

This case implies that $z=2a+1$.

Since $a+1=A(0,z)=\frac{1}{2}(z+1)$ and $z<b$, we have that $2a<b$ and
$A\cap [0,b]=[0,a]\cup\{b\}$ is a bi-arithmetic progression of difference $1$.

Since $|2(A\cap[a,n])|=3A(a,n)-3$, by applying Theorem \ref{3k-3} we have that
$A\cap [a,n]$ is either a bi-arithmetic progression, or $n-a+1\leqslant 2A(a,n)-1$,
or Freiman isomorphic to $K_6$ in (\ref{sixpoints}).

Notice that $n-a+1\leqslant 2A(a,n)-1=2A(z,n)+1=(n-z+1)+1$ implies that
$-a\leqslant -z+1$ and hence $2a+1=z\leqslant a+1$, which is absurd. So
we can assume that $A\cap [a,n]$ is either a bi-arithmetic progression
or Freiman isomorphic to $K_6$ in (\ref{sixpoints}).

\medskip

\quad {\bf Case 1.1.1}\quad $A\cap [a,n]$ is Freiman isomorphic to $K_6$
in (\ref{sixpoints}).

\medskip

Let $\varphi:K_6\mapsto A\cap [a,n]$ be the Freiman isomorphism.
Notice that $A(b+1,n-1)=3$.

Suppose that $b+1\not\in A$. Let $b'=\min A\cap[b+1,n]$.

If $a>1$, then there is an $x\in\{a-1,a-2\}$ such that $x+b'\not\in\{2a,a+b,2b\}$. Hence $x+b'$
is in the set in (\ref{twopieces}), which contradicts that the set is empty.

If $a=1$, then $z=3$, $b\geqslant 4$, and $n=12$ because $|A|=7$.
We can also assume that $b'=2b\geqslant 8$
because otherwise $0+b'$ is in the set in (\ref{twopieces}). Notice that
$a=1$ is a vertex of $\varphi(K_6)$ by part 1 of Proposition \ref{sixpoints2}.
Since $A\cap[a+1,b-1]=\emptyset$, $b$ is not a vertex by
part 2 of Proposition \ref{sixpoints2}. So $2b-a=2b-1$ is another
vertex of $\varphi(K_6)$. This contradicts the minimality of $b'$.

We can now assume that $b+1\in A$. If $b+1$ is a vertex of $\varphi(K_6)$,
then $c=\frac{1}{2}(a+b+1)$ is in $A$.
Clear, $a<c<b$, which contradicts $A\cap[a+1,b-1]=\emptyset$. So
$b+1$ is not a vertex in $\varphi(K_6)$. Hence $2b+2-a$ is in $A$ and is a vertex.
Since $2b-a$ and $2b+2-a$ are two vertices in $\varphi(K_6)$,
$2b+1-a$ is also in $A$. We can now conclude that
$A=[0,a]\cup\{b,b+1,2b-a,2b-a+1,2b-a+2\}$. So $(a-1)+(2b-a)=2b-1$ is in the
empty set in (\ref{twopieces}), which is absurd.

\medskip

\quad {\bf Case 1.1.2}\quad $A\cap [a,n]$ is a bi-arithmetic progression of difference $d$.

\medskip

Let $A\cap [a,n]=I_0\cup I_1$ be the bi-arithmetic progression
decomposition.

If $d=1$, then $A\cap[a,n]=\{a\}\cup [b,n]$ such that $n-b>b-a$.
Hence $A=[0,a]\cup[b,n]$ is a bi-arithmetic progression of difference $1$.

If $d=2$, then, without loss of generality,
$a\in I_1$ and $b\in I_0$ because $b-a\geqslant 3$.
Hence $I_1=\{a\}$.
But this contradicts the fact that $\gcd(A\cap[b,n])=1$.

If $d=3$, then $a,b$ should both be in $I_0$ or both be in $I_1$ because
$b-a\geqslant 3$ and $\gcd(A\cap[b,n]-b)=1$. Hence $b=a+3$ and $z=a+2$,
which implies that $a=1$ because $z=2a+1$.
Suppose, without loss of generality, $a\in I_0$. Let $c=\min I_1$.
If $c=b+2$, then $a-1+c$ is in the set in (\ref{twopieces}).
If $c=b+1$ or $c>b+3$, then $a-1+b+3$ is in
the set in (\ref{twopieces}). So both contradicts that the set in
(\ref{twopieces}) is empty.

If $d=4$, then $A(b+1,b+3)\leqslant 1$. If $a\not\equiv b\,(\mod 4)$,
then $b=a+3$ because $b-a\geqslant 3$ and $a+4=b+1\in A$. However, $b=a+3$
implies that $z=a+2=2a+1$ and $a=1$. So $A$ is a bi-arithmetic progression of difference $4$.
Hence we can assume that $a\equiv b\,(\mod 4)$. If $A(b+1,b+3)=0$ or
$b+1\in A$, let $x=b+4$. If $b+3\in A$, let $x=b+3$. Then $a-1+x$ is in
the empty set in (\ref{twopieces}), which is absurd.
Notice that $b+2\not\in A$ because otherwise
$\gcd(A\cap[b,n]-b)=2$.

If $d\geqslant 5$, then
\begin{eqnarray*}
\lefteqn{\frac{1}{2}(n-z+1)=A(z,n)=A(b,n)\leqslant\frac{2}{d}(n-b-1)+2}\\
& &\leqslant\frac{2}{5}(n-b-1)+2
=\frac{2}{5}(n-b+4)\leqslant\frac{2}{5}(n-z+3),
\end{eqnarray*}
which implies that $n-z\leqslant 7$.
Now $A(z,n)=\frac{1}{2}(n-z+1)$, $z<b$, and $d=5$ imply that
$n-z=7$, $d=5$, $b=z+1$, and $A\cap [b,n]=\{b,b+1,b+5,b+6\}$.
If $a\equiv b\,(\mod 5)$, then $a-1+b+5$ is in the empty set in
 (\ref{twopieces}). If $a\not\equiv b\,(\mod 5)$, then $b-a=4$,
 which implies that $a=2$ because $a+3=z=2a+1$.
 Hence $b+5=(a-2)+(b+5)=2b-1$ is in
 the empty set in (\ref{twopieces}). Both are absurd.

\medskip

{\bf Case 1.2}\quad $H\cap[0,a]\not=\emptyset$.

\medskip

Notice that $z\leqslant 2a$ in this case.
If $b>z+1$, then, by part 5 of Proposition \ref{addminim}, we have that
$|2(A\cap[0,b])|>3|A\cap[0,b]|-3$. Hence we can assume that $b=z+1$.

Since $A\cap [0,z]$ is anti-symmetric and $1\in A$, then $z-1\not\in A$. Hence $a<z-1$.
So $b-a\geqslant 3$.
If \[n-a+1\leqslant 2A(a,n)-1=2A(z,n)+1=(n-z+1)+1,\] then $z-1\leqslant a$.
Hence we can assume, by Theorem \ref{3k-3}, that $A\cap [a,n]$ is either a
bi-arithmetic progression or a Freiman isomorphism image of $K_6$ in (\ref{sixpoints}).

\medskip

\quad {\bf Case 1.2.1}\quad $A\cap [a,n]$ is Freiman isomorphic to $K_6$
in (\ref{sixpoints}).

\medskip

Let $\varphi:K_6\mapsto A\cap [a,n]$ be the Freiman isomorphism.

Since $A(z,n)=5=\frac{1}{2}(n-z+1)=\frac{1}{2}(n-b+2)$, we have that
$n-b=8$. Since $a$ is a vertex and $b$ is not a vertex of $\varphi(K_6)$,
we have that $2b-a$ is a vertex in $\varphi(K_6)$. Let $c$ be the third vertex
in $\varphi(K_6)$. If $2b-a=n$, then $c$ is between $n$ and $b$. Hence
$\frac{1}{2}(a+c)$ is in $A$ and is strictly between $a$ and $b$, which is impossible.
So we can assume that $2b-a<c=n$. Notice that
$\frac{1}{2}(n+2b-a)$ is in $A\cap[b,n]$. Clearly, $n-(2b-a)$ is even and $\leqslant 5$
because $n-b=8$ and $b-a\geqslant 3$.
If $n-(2b-a)=4$, then $A\cap[b,n]=\{b,b+2,b+4,b+6,b+8\}$, which contradicts
$\gcd(A\cap[b,n]-n)=1$. If $n-(2b-a)=2$, then $A\cap[b,n]=\{b,b+1,b+6,b+7,b+8\}$
and $a=b-6$. If $a-1\in A$, then $(a-1)+(b+6)$ is in the empty set in
(\ref{twopieces}). So we can assume that $a-1\not\in A$. Let
$a'=\max A\cap[0,a-1]$. If $a'+b+1\not= 2a$, then $a'+b+1$ is in the empty set in
(\ref{twopieces}). If $a'+b+1=2a$, then $a'+b+6$ is in the empty set in
(\ref{twopieces}). Both are absurd. This completes the proof
of Case 1.2.1.

\medskip

\quad {\bf Case 1.2.2}\quad $A\cap [a,n]$ is a bi-arithmetic progression of difference $d$.

\medskip

Let $A\cap[a,n]=I_0\cup I_1$ be the bi-arithmetic progression decomposition.

If $d=1$, then $A\cap[a,n]=\{a\}\cup [b,n]$ with $n-b<b-a$.
Let $A'=n-A$, $z'=n-z$, $b'=n-a$, and
$a'=n-b$. Then $A'\cap[0,z']=[0,a']$, $A'\cap[a',b']=[0,a']\cup\{b'\}$,
 and $z'$ is a balanced hole of $A'$. The same proof for Case 1.1
works for $A'$.

If $d=2$, then $\gcd(A\cap[b,n])=2$ because $b-a\geqslant 3$,
which contradicts the assumption that $\gcd(A\cap[b,n])=1$.

If $d=3$, then $a,b$ should both be in $I_0$ or both be in $I_1$ because
$b-a\geqslant 3$ and $\gcd(A\cap[b,n]-b)=1$. Hence $b=a+3$ and $z=a+2$.
Suppose, without loss of generality, $a\in I_0$. Let $c=\min I_1$ and
$a'=\max A\cap [0,a-1]$.

If $a'=a-1$, then $A(a',z)=2=\frac{1}{2}(z-a'+1)$. Hence $a'=0$
by the minimality of $z$, which implies that $H\cap[0,a]=\emptyset$, a contradiction
to the assumption of Case 1.2.

Thus we can assume that $a'<a-1$. If $a'\equiv a\,(\mod 3)$, then
$c+a'$ is in the empty set in (\ref{twopieces}). So we can assume that
$a'\not\equiv a\,(\mod 3)$. If $c=b+1$, then $c+a'$ is
in the empty set in (\ref{twopieces}). So we can assume that $c>b+1$.
If $b+3\not\in A$, then $c=b+2$ and
$A\cap [z,n]=\{b,b+2\}$ by the fact that $A(z,n)=\frac{1}{2}(n-z+1)$,
which contradicts the assumption that $\gcd(A\cap[b,n]-b)=1$.
So we can assume that $b+3\in A$. But now $a'+b+3$ is in the empty set in (\ref{twopieces}).

If $d=4$, then $b-a=4$ or $b-a=3$ because otherwise $A\cap[a,n]$ has a decomposition
$\{a\}\cup (A\cap[b,n])$, which contradicts $\gcd(A\cap[b,n]-b)=1$.
Let $a'=\max(A\cap[0,a-1])$ and $c=\min\{x\in A\cap[b+1,n]:x\not\equiv b\,(\mod 4)\}$.

Suppose that $a'=a-1$. If $b-a=4$, then $c\not=b+2$ and $b+4\in A$
by the fact that $A(b-1,n)=\frac{1}{2}(n-b)$. If $c\not=b+3$, then $a'+b+4$
is in the empty set in (\ref{twopieces}). If $c=b+3$, then
$a'+c$ is in the empty set in (\ref{twopieces}).
If $b-a=3$, then $z-a=2$ and $A\cap[a',z]=\{a',a\}$.
Hence $A(a',z)=\frac{1}{2}(z-a'+1)$, which implies
that $a'=0,a=1$ by the minimality of $z$. But this contradicts the assumption that
$[0,a]\cap H\not=\emptyset$.

So we can now assume that $a'<a-1$. If $b+1\in A$, then $a'+b+1$
is in the empty set in (\ref{twopieces}) unless
$a'+b+1=2a$. If $a'+b+1=2a$, then $2a\not\in (b+A\cap[0,b])$,
which leads to a contradiction to (\ref{sub3k-3}).
So we can assume that $b+1\not\in A$, which implies that $b\not=a+3$
because otherwise $a+4=b+1\in A$.
So we have that $b=a+4$. Since $A(z,n)=\frac{1}{2}(n-z+1)$, we
have that $b+3\in A$ and $A\cap [a,n]$ has the decomposition
$\{a,a+4,\ldots,n\}\cup\{b+3,b+7,\ldots,n-1\}$ with $n\geqslant b+4$
by the fact that $A(z,n)=\frac{1}{2}(n-b)$. If $a'\equiv a\,(\mod 4)$, then
$a'+b+3$ is in the empty set in (\ref{twopieces}).
If $a'\not\equiv a\,(\mod 4)$, then $a'+b+4$
is in the empty set in (\ref{twopieces}).

If $d\geqslant 5$, then
\begin{eqnarray*}
\lefteqn{\frac{1}{2}(n-z+1)=A(z,n)=A(b,n)}\\
& &\leqslant\frac{2}{d}(n-b-1)+2\leqslant\frac{2}{5}(n-b-1)+2\\
& &=\frac{2}{5}(n-b+4)\leqslant\frac{2}{5}(n-z+3),
\end{eqnarray*}
which implies that $n-z\leqslant 7$. Since $A(b-1,n)=\frac{1}{2}(n-b)$,
we have that $n-z=7$, $d=5$, and $A\cap [b,n]=\{b,b+1,b+5,b+6\}$.
Notice that $b-a$ is $5$ or $4$. Let $a'=\max A\cap[0,a-1]$.

If $a'<a-1$, then $a'+b+1$ is in the empty set
 in (\ref{twopieces}) unless
$a'+b+1=2a$. If $a'+b+1=2a$, then $2a\not\in (b+A\cap[0,b])$, which leads to a
contradiction to (\ref{sub3k-3}).

So we can assume that $a'=a-1$.
If $a-b=5$, then $a'+b+5$ is the empty set in (\ref{twopieces}).
If $a-b=4$, then let $a''=\max\{x\in A\cap[0,a-1]:x\not\equiv b,b+1\}$.
Notice that
$a''$ exists because otherwise $A$ is a subset of a bi-arithmetic progression
of difference $5$, which leads to a contradiction to the assumption that
$H\cap [0,a]\not=\emptyset$ and $A\cap[0,z]$ is left dense.
If $a''\equiv b+2$ or $b+3\,(\mod 5)$, then $a''+b+1$ is in the empty set in
(\ref{twopieces}). If $a''\equiv b+4\,(\mod 5)$, then
$a''+b+5$ is the empty set in (\ref{twopieces}).

This completes the proof of Case 1.2.2 as well as Case 1.

\medskip

\noindent {\bf Case 2}\quad $1\not\in A$.

\medskip

By symmetry, we can also assume that $n-1\not\in A$.
We want to show that $A$ is in the form of $T_{k,n}$ or $n-T_{k,n}$
for some $k>0$ defined in (\ref{tk}) or $A$ is in the form of $S_{u,n}$ defined in (\ref{su}).
Let $E$ be the set of all even numbers and
\[a=\max\{x\in [0,n]:A\cap[0,x]=E\cap [0,x]\}.\] Notice that $0<a<n$. Notice also that
if $a\in A$, then $a+1\in A$, and if $a\not\in A$,
then $a+1\not\in A$ by the maximality of $a$. If $a\in A$, we show that
$A$ is in the form of $T_{k,n}$ for $k=a/2$. If $a\not\in A$, then we show that
$a>1$ implies $|2A|>3|A|-3$ and $a=1$ implies that $A$ is either in the form of
$n-T_{k,n}$ or in the form of $S_{u,n}$.

\medskip

{\bf Case 2.1}\quad $a\in A$.

\medskip

Let $A'=A\cap [a,n]$. Then $a,a+1\in A$. Notice that
\begin{eqnarray*}
\lefteqn{3|A|-3=|2A|}\\
& &\geqslant |2(A\cap[0,a])|+|a+1+A\cap[0,a-2]|+|2(A\cap[a,n])|-1\\
& &=2A(0,a)-1+A(0,a-2)+|2(A\cap[a,n])|-1\\
& &=3A(0,a-2)+|2(A\cap[a,n])|.
\end{eqnarray*}
Hence $|2(A\cap[a,n])|\leqslant 3A(a,n)-3$. Notice also that
\[A(a,n)=|A|-A(0,a-1)=\frac{1}{2}(n+2)-\frac{1}{2}a=\frac{1}{2}(n-a+2).\]
So $|2(A\cap[a,n])|=3A(a,n)-3$ by part 3 of Proposition \ref{sub3k-3prop}.
Let $n'=n-a$. Now $A'=A\cap [a,n]-a$ in $[0,n']$ satisfies all conditions for Case 1.
So either $A'$ is in the form of $T_{0,n'}$ or $A'$ is a
bi-arithmetic progression of difference $1$ or $4$. Since $A$ is assumed not to be
a bi-arithmetic progression, we regard it as a contradiction when $A$ is forced to
be a bi-arithmetic progression. However, $A'$ is not $A$. So it is possible that
$A'$ is a bi-arithmetic progression of difference $1$ or $4$.

Suppose that $A'$ is a bi-arithmetic progression of difference $1$.
Let $A\cap [a,n]=[a,x]\cup [y,n]$
be the decomposition. Since $n-1\not\in A$,
we have that $y=n$. Since $|A\cap [a,n]|=|A|-A(0,a-1)
=\frac{1}{2}(n+2)-\frac{1}{2}a=\frac{1}{2}(n-a+1)$,
we have that $A'$ is in the form of $T_{0,n-a}$.
Hence $A$ is in the form of $T_{\frac{a}{2},n}$.

Suppose that $A'$ is a bi-arithmetic progression of difference $4$.
If $a+5\not\in A$, then
$A'=\{0,1,4\}$ by the fact that $|A'|=\frac{1}{2}(n'+2)$.
But $\{0,1,4\}$ is a form of $T_{0,4}$. Hence $A$ is in the form of
$T_{\frac{a}{2},a+4}$.
So we can assume that $a+5\in A$. Now we have the following contradiction:
\begin{eqnarray*}
\lefteqn{3|A|-3=|2A|\geqslant |2(A\cap[0,a])|+|a+1+A\cap[0,a-2]|}\\
& &\quad +|\{(a+5)+(a-2)\}|+|2(A\cap[a,n])|-1\\
& &\geqslant 3A(0,a-2)+1+3A(a,n)-3=3|A|-2>3|A|-3
\end{eqnarray*}
Therefore, we can conclude that $A'$ is in the form of $T_{0,n'}$ and
hence $A$ is in the form of $T_{k,n}$ for $k=a/2$.

\medskip

{\bf Case 2.2}\quad $a\not\in A$.

\medskip

So $a+1\not\in A$. Notice that $A(0,a)=\frac{1}{2}(a+1)$. Thus $a$
is a balanced hole. Notice also that $A(a,n-2)=\frac{1}{2}(n-a-1)$ because
$n\in A$ and $n-1\not\in A$. Let
\begin{equation}\label{definitionzu}
u=\min\left\{x\in [a,n-2]:A(a,x)\geqslant\frac{1}{2}(x-a+1)\right\}.
\end{equation}
Notice that $u>a+2$ because $a,a+1\not\in A$. Notice also that
$u,u-1\in A$, $A(a,u)=\frac{1}{2}(u-a+1)$, and
$A(x,u)>\frac{1}{2}(u-x+1)$ for every $x\in [a+1,u]$
by the minimality of $u$. So $A\cap[a,u]$ is right dense in $[a,u]$ and
\begin{eqnarray*}
\lefteqn{A(u,n)=|A|-A(0,a-1)-A(a,u-1)}\\
& &=\frac{1}{2}(n+2)-\frac{1}{2}(a+1)-\frac{1}{2}(u-a-1)=\frac{1}{2}(n-u+2).
\end{eqnarray*}
By part 3 of Proposition \ref{sub3k-3prop}, we can assume that
$|2(A\cap[0,u])|\geqslant 3A(0,u)-3$.
Notice that if $n'=u$, $z'=u-a$, and $A'=u-A\cap [0,u]$, then $A'\cap[0,z']$
is left dense in $[0,z']$. Applying the proof of Case 1 to $A'$,
we have that if $|2(A\cap[0,u])|=3A(0,u)-3$, then
either $A\cap [0,u]$ is in the form of $u-T_{0,u}$,
which is possible only when $a=1$,
or $A\cap [0,u]$ is a bi-arithmetic progression of difference $4$.

\medskip

\quad {\bf Case 2.2.1}\quad $a>1$.

\medskip

In this case we derive a contradiction by showing $|2A|>3|A|-3$.

Since $0,2\in A$ and $1\not\in A$, $A$ as well as $A\cap [0,u]$
can be neither a bi-arithmetic progression of
difference $1$ nor a bi-arithmetic progression of difference $4$.
So we want to derive a contradiction by show that $|2A|>3|A|-3$.
Notice that since $0,2\in A$ and $1\not\in A$, $A\cap[0,u]$ cannot be
in the form of $u-T_{0,u}$. Hence by applying the proof of Case 1 to $u-(A\cap[0,u])$,
we have that
\[|2(A\cap[0,u])|>3A(0,u)-3.\]

If $\gcd(A\cap[u,n]-u)>1$, then
\begin{eqnarray*}
\lefteqn{|2A|\geqslant|2(A\cap [0,u])|+|2(A\cap[u,n])|-1
+|u-1+A\cap[u+2,n]|}\\
& &>3A(0,u)-3+2A(u,n)-2+A(u+1,n)=3|A|-3.
\end{eqnarray*}
Hence we can assume that $\gcd(A\cap[u,n]-u)=1$, which implies that
$|2(A\cap[u,n])|\geqslant 3A(u,n)-3$.

Let $v=\min A\cap[u+1,n]$.

If $v>u+1$, then $u-1+v$ is in the set
\begin{equation}\label{twopieces2}
(2A)\smallsetminus(2(A\cap[0,u]))\cup(2(A\cap[u,n])).
\end{equation}
Hence
\begin{eqnarray*}
\lefteqn{|2A|\geqslant |2(A\cap[0,u])|+|2(A\cap[u,n])|}\\
& &>3A(0,u)-3+3A(u,n)-3=3|A|-3.
\end{eqnarray*}
Thus we can assume that $v=u+1$.

Recall that in the proof of $\gcd(A\cap[b,n]-b)=1$ at the beginning of
Case 1 but before case 1.1, we showed that if $b<n$ and $\gcd(A\cap [b,n]-b)>1$,
then $|2A|>3|A|-3$ without counting the possible element $z=b-1$.
Recall also that $\gcd(A\cap[0,a])=2$. So
we can use the same argument to $A'=u-(A\cap [0,u])$ with $b'=u-a+1$
to show that
\[|2(A\cap[0,u])\cup\{u+a\}|=|(2A')\cup\{b'-1\}|>3|A'|-3+1=3A(0,u)-2.\] Hence
\begin{eqnarray*}
\lefteqn{|2A|\geqslant
|2(A\cap[0,u])\cup\{a-1+v\}|+|2(A\cap [u,n])|-1}\\
& &>3A(0,u)-2+3A(u,n)-4=3|A|-3.
\end{eqnarray*}

\medskip

\quad {\bf Case 2.2.2}\quad $a=1$.

\medskip

In this case we show that $A$ is either in the form of
$n-T_{k,n}$ for some $k>0$ or in the form
of $S_{u,n}$ for some $u\in [4,n-6]$.

Notice that $A\cap [1,u]$ is a right dense set in $[1,u]$
and $|2(A\cap[0,u])|\geqslant 3A(0,u)-3$
by part 3 of Proposition \ref{sub3k-3prop}.
If $|2(A\cap[0,u])|=3A(0,u)-3$, then $A\cap[1,u]$ is
a right dense, additively minimal, anti-symmetric set in $[1,u]$.
If $A\cap [u,n]$ is an arithmetic progression of
difference $2$, then $A$ is in the form of $n-T_{k,n}$ for some $k>0$.
So we can assume $\gcd(A\cap [u,n]-u)=1$
if $|2(A\cap[0,u])|=3A(0,u)-3$. But we haven't eliminate the possibility that
$|2(A\cap[0,u])|>3A(0,u)-3$ can happen.

\medskip

\quad\quad {\bf Case 2.2.2.1}\quad $u+1\in A$.

\medskip

If $|2(A\cap[0,u])|>3A(0,u)-3$, then
\begin{eqnarray*}
\lefteqn{3|A|-3=|2A|\geqslant |2(A\cap[0,u])|+|2(A\cap[u,n])|-1}\\
& &\geqslant 3A(0,u)-2+3A(u,n)-3-1=3|A|-3.
\end{eqnarray*}
Hence $|2(A\cap[u,n])|=3A(u,n)-3$. By applying the proof of Case 1 to the
set $A'=A\cap[u,n]-u$ and $n'=n-u$, we can conclude that
$A'$ is either in the form of $T_{0,n'}$ or a bi-arithmetic progression of difference $1$ or $4$.
We now want to show that $|2A|>3|A|-3$ by identifying one element in the set in
(\ref{twopieces2}), which implies that $|2A|>3A(0,u)-3+3A(u,n)-3-1+1=3|A|-3$.

If $A'$ is in the form of $T_{0,n'}$, then $u-1+n$ is in the set in (\ref{twopieces2}).
If $A'$ is a bi-arithmetic progression of $1$, then
$A\cap[u,n]=[u,x]\cup\{n\}$ because $u-1\not\in A$. So again $A'$ is in the form of
$T_{0,n'}$.
If $A'$ is a bi-arithmetic progression of difference $4$, then $u-1+u+4$
is in the set in (\ref{twopieces2}).

Thus we can assume that $|2(A\cap[0,u])|=3A(0,u)-3$. So
$A'=u-(A\cap[0,u])$ is in the form of $T_{0,n'}$ for $n'=u$ or
$A'$ is a bi-arithmetic progression of difference $1$ or $4$. Notice that
if $A'$ is a bi-arithmetic progression of difference $1$, then $A'$ is in the form of $T_{0,n'}$
because $1\not\in A$. And
if $A'$ is a bi-arithmetic progression of difference $4$, then $A'=\{0,1,4\}$
which is also in the form of $T_{0,n'}$,
because $A'$ is left dense in $[0,u-1]$ and $A'(0,3)=\frac{1}{2}(n'+1)$.
As a consequence we have that $u+1$ is in the set in (\ref{twopieces2}).

If $|2(A\cap[u,n])|>3A(u,n)-3$, then $|2A|>3A(0,u)-3+3A(u,n)-3=3|A|-3$.
Hence we can now assume that $|2(A\cap[u,n])|=3A(u,n)-3$.

Now we have assumed that $|2(A\cap[0,u])|=3A(0,u)-3$, $|2(A\cap[u,n])|=3A(u,n)-3$,
$\{u-1,u,u+1\}\subseteq A$, and $A\cap[0,u]$ is in the form of $u-T_{0,u}$.
By applying the proof of Case 1 to $A'=(A\cap [u,n])-u$, we have that
$(A\cap[u,n])-u$ is either in the form of $T_{0,n-u}$ or
a bi-arithmetic progression of difference $4$. Notice that
if $A\cap [u,n]$ is a bi-arithmetic progression of difference $1$,
then $A\cap [u,n]-u$ is in the form of $T_{0,n-u}$.
We now want to show that $|2A|>3|A|-3$ by identifying two elements in the set in
(\ref{twopieces2}), which implies that $|2A|\geqslant 3A(0,u)-3+3A(u,n)-3-1+2=3|A|-2$.

If $(A\cap[u,n])-u$ is in the form of $T_{0,n-u}$, then
$u+1,u-1+n$ are in the set in (\ref{twopieces2}).
If $A\cap[u,n]$ is a bi-arithmetic progression of difference $4$, then
$0+u+1,u-1+u+4$ are in the set in (\ref{twopieces2}).

\medskip

\quad\quad {\bf Case 2.2.2.2}\quad $u+1\not\in A$.

\medskip

Let $v=\min A\cap[u+1,n]$. Notice that $v+u-1$ is in the set in (\ref{twopieces2}).
If $|2(A\cap[0,u])|>3A(0,u)-3$, then $|2A|>3A(0,u)-3+3A(u,n)-3=3|A|-3$.
Hence we can assume that $|2(A\cap[0,u])|=3A(0,u)-3$. By applying the proof of Case 1,
we have that $u-(A\cap[0,u])$ is in the form of $T_{0,u}$.
If $\gcd(A\cap [u,n]-u)=d>1$, then $d=2$ and $A\cap [u,n]$ is an arithmetic progression
of difference $2$. So $A$ is in the form of $n-T_{k,n}$ for $k=(n-u)/2$. Hence
we can assume that $\gcd(A\cap [u,n]-u)=1$. If $|2(A\cap[u,n])|>3A(u,n)-3$,
then $|2A|>3A(0,u)-3+1+3A(u,n)-3-1=3|A|-3$. Hence we can assume that
$|2(A\cap[u,n])|=3A(u,n)-3$.

Suppose that $v=u+2$. We want to show that $A$ is in the form of $S_{u,n}$.

Let $a'=\max\{x\in [u,n]:A\cap[u,x]=(u+E)\cap[u,x]\}$. $a'$ is well defined because
$\gcd(A\cap[u,n]-u)=1$. Notice that $a'\geqslant u+2$.

If $a'\not\in A$, then $a'>u+2$. By the proof of Case 2.2.1 we have that
$|2(A\cap[u,n])|>3A(u,n)-3$, a contradiction. So we can assume that $a'\in A$,
which implies that $a'+1\in A$ by the maximality of $a'$.
By applying the proof of Case 2.1
to the set $A'=A\cap[u,n]-u$, we have that $A'$ is in the form of $T_{k,n-u}$.
 Hence
\[A=\{0\}\cup C\cup\{u,u+2,\ldots,u+2k\}\cup D\cup\{n\}\]
where $C$ is right dense, anti-symmetric, and additively minimal in $[1,u]$, and
$D$ is left dense, anti-symmetric, and additively minimal in $[u+2k,n-1]$.
Thus
\begin{eqnarray*}
\lefteqn{|2A|\geqslant |2(A\cap[0,u])|+|[2u+1,2u+4k-1]|+|2(A\cap[u+2k,n])|}\\
& &\geqslant 3A(0,u)-3+4k-1+3A(u+2k,n)-3\\
& &=3A(0,u)-3+4A(u+2,u+2k-2)+3+3A(u+2k,n)-3\\
& &=3|A|-3+A(u+2,u+2k-2)>3|A|-3
\end{eqnarray*}
unless $k=1$. If $k=1$, then $A$ is in the form of $S_{u,n}$ in (\ref{su}).

Now we assume that $v>u+2$. We show that $|2A|>3|A|-3$ by identifying two elements
in the set in (\ref{twopieces2}).

If $u-2,u-3\not\in A$, then $u=4$ and $A\cap[0,u]=\{0,3,4\}$ by the minimality of $u$.
Since $A$ is not a bi-arithmetic progression, we can define
\begin{eqnarray*}
\lefteqn{c=\min\{x\in [u,n]:A\cap[0,x]}\\
& &\mbox{ is not a subset of a bi-arithmetic progression of
difference }\,4\}.
\end{eqnarray*}
Since $3,4\in A$, we have that
$c$ is either congruent to $5$ or congruent to $6$ modulo $4$.

Suppose that $c=v$. Recall that $v>6$.
So we can assume that $v\geqslant 9$. Then
$v+3,v+0$ are in the set in (\ref{twopieces2}).

Suppose that $c>v$. Then $v$ is congruent to $3$ or $4$ modulo $4$.
If $c\equiv 5\,(\mod 4)$, then $c+0,v+3$ are in the set in (\ref{twopieces2}).
If $c\equiv 6\,(\mod 4)$, then
$c+3,v+3$ are in the set in (\ref{twopieces2}).

So we can assume that $A(u-3,u)\geqslant 3$.
If $u-2\in A$, then $v+u-1,v+u-2$ are in the set in (\ref{twopieces2}).
So we can assume that $u-2\not\in A$ and $u-3\in A$.

If $v\geqslant u+4$, then $v+u-1,v+u-3$ are in the set in (\ref{twopieces2}).
So we can assume that $v=u+3$. Let $v'=\min A\cap[v+1,n]$.
If $v'=v+1$, then $v+u-1,v'+u-3$ are in the set in (\ref{twopieces2}).
If $v'>v+1$, then $v+u-1,v'+u-1$ are in the set in (\ref{twopieces2})
unless $v'+u-1=2v$. But if
$v'+u-1=2v$, then $v+u-1,v'+u-3$ are in the set in (\ref{twopieces2}).

This completes the proof of Theorem \ref{main1}.

\section{Questions}

We end this paper by asking few questions. The first question is related to Theorem \ref{BGthm}.

\begin{question}
Let $A$ and $B$ be finite nonempty sets of integers with \[\max B-\min B\leqslant
\max A-\min A\leqslant |A|+|B|-3\]
\[\mbox{and }\,|A+B|= |A|+2|B|-2-\delta(A,B).\]
What should be the structure of $A$, $B$, and $A+B$, which generalizes both Theorem \ref{BGthm} and
Theorem \ref{main1}?
\end{question}

The second question is one step further than Theorem \ref{main1}.

\begin{question}
Let $A$ be a finite nonempty set of integers such that $|A|\geqslant 11$ and
\[|2A|=3|A|-2.\]
What should be the detailed structure of $A$ and $2A$, which generalizes Theorem \ref{main1}?
\end{question}

Notice that we define a left dense and anti-symmetric set $B$ to be additively minimal in $[u,v]$
by using a property of the sumset $2(B\cup\{v+1\})$. The following question is vague on purpose.

\begin{question}
What should be a nice direct characterization of a left dense, anti-symmetric, additively minimal set
$B$ without mentioning the sumset.
\end{question}

\end{document}